\documentclass[10pt]{article}
\usepackage[pagewise]{lineno}
\usepackage{bbm}
\usepackage{amsmath}
\usepackage{tabularx}
\usepackage{enumitem}
\usepackage{amssymb} 
\usepackage{mathrsfs}
\usepackage{caption}
\usepackage{graphicx}
\usepackage{authblk}
\usepackage[hypertexnames=false,colorlinks=true,linkcolor=blue,citecolor=blue]{hyperref}
\usepackage[numbers,comma,square,sort&compress]{natbib}
\usepackage[letterpaper,text={6.5in,9in},centering]{geometry}

\graphicspath{{eps/}{pdf/}}

\makeatletter\@addtoreset{equation}{section}\makeatother
\renewcommand{\theequation}{\arabic{section}.\arabic{equation}}

\newtheorem{Lemma}{Lemma}[section]
\newtheorem{Theorem}{Theorem}

\newtheorem{Remark}[Lemma]{Remark}

\newenvironment{Proof}%
 {\begin{trivlist} \item[]{\bf Proof. }}%
 {\hspace*{\fill}$\rule{.4\baselineskip}{.4\baselineskip}$\end{trivlist}}

\newenvironment{Acknowledgment}%
 {\begin{trivlist}\item[]\textbf{Acknowledgments.}}{\end{trivlist}}
\def\Re{\mathop{\mathrm{Re}}}

\usepackage{tikz}

\begin{document}
\title{Existence of localised radial patterns in a model for dryland vegetation}
\author[1,2]{Dan J. Hill}

\affil[1]{\small Department of Mathematics, University of Surrey, Guildford, GU2 7XH, UK}
\affil[2]{\small Fachrichtung Mathematik, Universit\"at des Saarlandes, Postfach 151150, 66041 Saarbr\"ucken, Germany}
\date{\today}
\maketitle

\begin{abstract}
\noindent
Localised radial patterns have been observed in the vegetation of semi-arid ecosystems, often as localised patches of vegetation or in the form of `fairy circles'. We consider stationary localised radial solutions to a  reduced model for dryland vegetation on flat terrain. By considering certain prototypical pattern-forming systems, we prove the existence of three classes of localised radial patterns bifurcating from a Turing instability. We also present evidence for the existence of localised gap solutions close to a homogeneous instability. Additionally, we numerically solve the vegetation model and use continuation methods to study the bifurcation structure and radial stability of localised radial spots and gaps. We conclude by investigating the effect of varying certain parameter values on the existence and stability of these localised radial patterns.
\end{abstract}

\section{Introduction}
\label{s:intro}

\begin{figure}[t!]
    \centering
    \includegraphics[width=\textwidth]{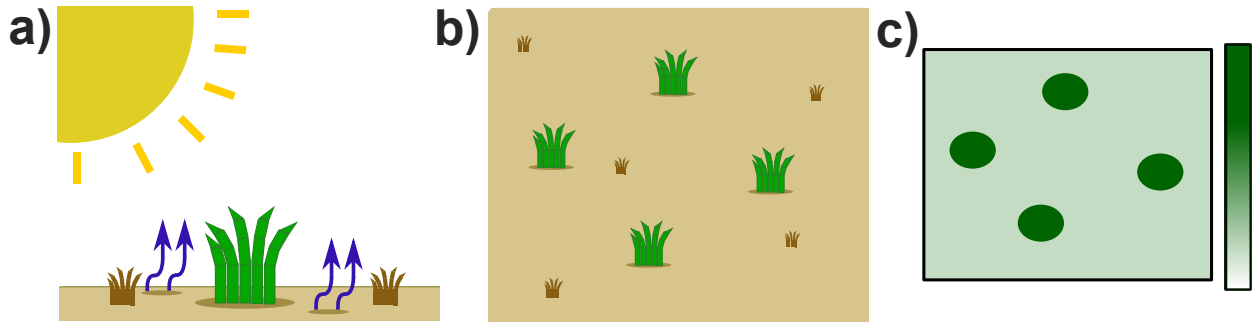}
    \caption{a) Due to evaporation of the soil water (illustrated by the blue arrows) and insufficient resources to sustain uniform plant growth, semi-arid ecosystems can exhibit localised patches of vegetation. The soil water is also localised, subsisting in the shadow of the vegetation. b) For sufficiently large domains, these patches can also form part of a larger pattern (such as spots, stripes or hexagons, for example). c) By measuring the vegetation and soil water densities, vegetation patterns can be described by continuum models. Here, we present a schematic density plot of the vegetation in panel b), where the darker green indicates a higher intensity of vegetation.}
    \label{fig:cartoon}
\end{figure}

The self-organisation of vegetation in semi-arid ecosystems remains one of the most visually striking examples of pattern formation in nature. 
Non-uniform distributions of vegetation are thought to result from symmetry-breaking instabilities induced by a lack of resources, such as a scarcity of water, combined with a positive feedback between local vegetation growth and water transport towards the growing vegetation; see \cite{MERON2016Pattern} for a detailed review of the mechanisms involved in the formation of vegetation patterns, and see Figure \ref{fig:cartoon} for a schematic of how vegetation patterns form and are measured. Semi-arid ecosystems are defined to have an annual rainfall of $250$-$500$ mm; these ecosystems are not rare, with $90\%$ of Israel defined to be semi-arid or arid, and $30\%$ of the Earth's land surface estimated to be semi-arid \cite{vanderstelt2013}. With the advent of aerial \cite{MacFadyen1950a,MacFadyen1950b} and satellite imagery, and an increasing public awareness of ecological issues such as desertification, the study of vegetation patterns has never been more accessible or relevant; this is reflected in the wealth of new research being completed in this area  (see \cite{Berrios-Caro2020,Gandhi2020,Jaibi2020,Tlidi2008,Siero2020} for recent work in mathematical modelling, \cite{YOUSSEF2012Wind,McGrath2012Microtopography} for experimental and numerical simulations, and \cite{Rasanen2020Peatland,Deblauwe2011Modulation,Deblauwe2012Migration} for remote sensing from satellite imagery).

While sloped terrains often produce slowly migrating stripe patterns \cite{Samuelson2019,Siero2015,vanderstelt2013}, flat terrains are known to produce a variety of different steady structures \cite{Borgogno2009,Gowda2014,Gowda2016,Rietkerk2002}. These patterns include localised patches, where sparse regions of vegetation are surrounded by bare soil \cite{Meron2007}, and `fairy circles', where a sparse region of bare soil is surrounded by uniform vegetation (as seen in the Namib desert \cite[Figure 1.]{Tlidi2008}). A single patch or fairy circle is locally axisymmetric and spatially localised; hence, we restrict our focus to localised radial patterns in order to better understand these structures.

There has been a number of analytical and numerical studies of localised vegetation patterns; see \cite{Meron2007,Tlidi2008,Berrios-Caro2020,Bordeu2016} for nonlocal mean-field models, \cite{Berrios-Caro2020,Bordeu2016,Lejeune2002} for non-variational PDE models, and \cite{Escaff2015} for a variational nonlocal Nagumo-type model. Notably, Ja\"ibi et al. \cite{Jaibi2020} rigorously analysed a 2-component reaction-diffusion model in order to prove the existence of localised travelling fronts, as well as stationary spot and gap patterns. This work focused on one-dimensional patterns and utilised singular perturbation theory in order to complete the proof; it may be possible to extend these results for the radial case, however the associated equations will explicitly depend on the radial coordinate, which may introduce some extra technical challenges. Instead, we will take a weakly nonlinear approximation of the reaction-diffusion model introduced by von Hardenberg et al. in \cite{vonHardenberg2001}, which we term the `von Hardenberg' model, and draw connections to prototypical pattern-forming models, which we will discuss presently.

\begin{figure}[t!]
    \centering
    \includegraphics[width=\textwidth]{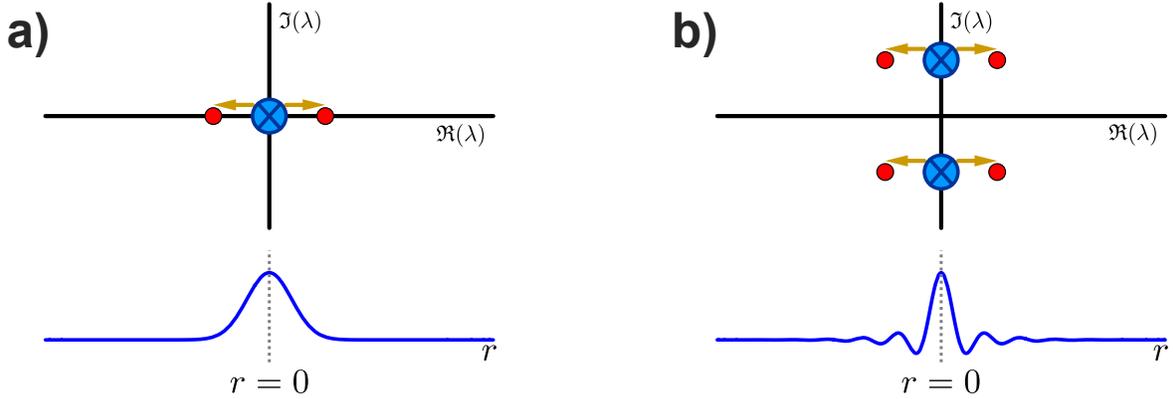}
    \caption{a) For localised spikes, the top picture shows a pair of real eigenvalues bifurcating from the origin (blue circle with cross) onto the real line (red circles). This results in a localised solution with monotonic tails, as shown in the bottom picture. b) For localised spots, the top picture shows four purely imaginary eigenvalues bifurcating from $\pm \textnormal{i} k$ (blue circle with cross), for some $k\in\mathbb{R}$, into the complex plane (red circles). This results in a localised solution with spatially-oscillatory tails with wave number $k$, as shown in the bottom picture.}
    \label{fig:Pulses}
\end{figure}

Within the study of localised states, we wish to highlight an important distinction between \textbf{localised spikes} and \textbf{localised spots}; here we are using the terminology set out in \cite{Champneys2021}, related to the spatial eigenvalues of the surrounding homogeneous state at the bifurcation point. We note that, for radial problems, we typically examine the spatial eigenvalues in the limit as $r\to\infty$. A localised spike is commonly the result of spatial eigenvalues bifurcating from the origin onto the real line; these solutions do not have an associated wave number, and so they have monotonic decay in their tails; see Figure \ref{fig:Pulses} a). In contrast, a localised spot is often the result of purely imaginary eigenvalues bifurcating off of the imaginary axis; such solutions have an associated wave number, and so they decay with spatially-oscillating tails; see Figure \ref{fig:Pulses} b). Localised patterns associated with this bifurcation often possess much more complicated bifurcation curves than localised spikes and are known to sometimes undergo homoclinic snaking, where infinite families of localised patterns are connected via a snaking process, as seen in \cite{Woods1999}. For more information about homoclinic snaking, see \cite{burke2007homoclinic} and references therein.

In the context of vegetation patterns, localised spikes are typically found far from any symmetry-breaking instability, while localised spots are found close to such an instability (see \cite{Champneys2021} for a review of such structures). This instability is often called a `Turing' bifurcation after Turing's seminal work on pattern formation \cite{turing1952chemical}. Since we will be studying a reaction-diffusion model, such that the instabilities of the system are diffusion-driven, we will stick to this naming convention; we will also refer to the spike-inducing bifurcation, where real eigenvalues are bifurcating from the origin, as a `homogeneous' bifurcation. A notable difference between localised spikes and spots can be seen in how they interact; see \cite{Berrios-Caro2020,Tlidi2008,Tlidi2020,Bordeu2016}. Berrios-Caro et al. \cite{Berrios-Caro2020} found that, for a nonlocal reaction-diffusion equation, localised spikes can be either purely attractive or purely repulsive; in contrast, the interaction between localised spots alternates between attractive and repulsive depending the distance separating the patterns, allowing for the existence of bound states. Similar results can be found in \cite{Tlidi2020}, where localised spots with monotonic tails are seen to be repulsive while localised gaps with oscillatory tails were found to exhibit bound states.  Throughout this work, we will attempt to distinguish between localised spikes and localised spots, identifying their bifurcation points and regions of existence for the von Hardenberg model. However, we first require some background theory regarding the existence of localised radial patterns, which we will now discuss in the context of some prototypical models.

\paragraph{Localised Radial Patterns} 
One of the biggest issues in the analytical study of localised planar patterns is when localisation occurs in multiple directions. There has been some progress regarding two-dimensional patterns, such as \cite{buffoni2018variational,HillPrepPatch,sandstede2020spiral}, but the majority of analytical work focuses on localised radial patterns. In this case the localisation is restricted to just the radial direction, effectively reducing the problem to one independent variables. Localised radial patterns are seen as the first step to better understanding localised planar patterns \cite{knobloch2008spatially}; however, radial problems tend to have explicit dependence on the spatial coordinate $r$, which can cause significant problems in the analysis. 

There are a number of papers that investigate the existence of localised radial patterns in various systems, such as \cite{vanHeijster2011,Hill2020Localised,mcquighan2014oscillons,mccalla2013spots,lloyd2009localized}. However, all of these examples utilise the radial invariant manifold and normal form theory introduced by Scheel \cite{scheel2003radially}, and so we quickly summarise the relevant results from \cite{scheel2003radially} and their applications in two prototypical examples.

\underline{\textit{General results for localised radial patterns:}} We begin with a general stationary axisymmetric reaction-diffusion system,
\begin{align}
    \mathbf{D}\Delta_{r} \mathbf{U} + \mathbf{F}(\mathbf{U};\mu) = 0,\nonumber
\end{align}
where $\mathbf{D}$ is a diffusion matrix, $\Delta_{r}:=\partial_{rr} + \frac{1}{r}\partial_{r}$ is the radial Laplacian operator, $\mathbf{U}\in\mathbb{R}^{2}$ is the dependent variable, $\mu\in\mathbb{R}^{p}$ is a parameter, and $\mathbf{F}$ is a nonlinear function where $\mathbf{F}(0;0)=0$. We can also express this as a first-order ODE system 
\begin{align}
    \mathbf{u}_{r} &= \mathcal{A}(r)\mathbf{u} + \mathcal{F}(\mathbf{u};\mu), \label{Scheel:Syst}\\
    \intertext{where $\mathbf{u}=(\mathbf{U},\partial_{r}\mathbf{U})^{\intercal}\in\mathbb{R}^{4}$,}
    \mathcal{A}(r) &= \begin{pmatrix}
    0 & 1_{2} \\ -\mathbf{D}^{-1}\partial_{\mathbf{U}}\mathbf{F}(0;0) & -\frac{1}{r}\,1_{2}
    \end{pmatrix}, \qquad \mathcal{F}(\mathbf{u};\mu) = \begin{pmatrix}
    0 \\ -\mathbf{D}^{-1}\left(\mathbf{F}(\mathbf{U};\mu) - \partial_{\mathbf{U}}\mathbf{F}(0;0)\right)
    \end{pmatrix},\nonumber
\end{align}
where $1_{2}$ denotes the $2$-dimensional identity matrix. The key idea is the following; a localised radial pattern remains bounded close to the origin, which we will call the `core' region, while exhibiting exponential decay sufficiently far from the origin, which we will call the `far-field' region. The existence proof follows from a singular perturbation argument, and involves local invariant manifold theory. We first construct the set of all possible small-amplitude solutions that remain bounded in the core region. We call this set of solutions the core manifold, denoted by $\mathcal{W}^{cu}_{-}(\mu)$, which we construct over the finite sub-interval $r\in[0,r_{0}]$. Likewise, we then construct the set of all small-amplitude solutions that exhibit exponential decay in the far-field region. We call this set of solutions the far-field manifold, denoted by $\mathcal{W}^{s}_{+}(\mu)$, which is defined for sufficiently large values of $r$. We note that we are borrowing the notation for standard `centre-unstable' and `stable' manifolds to denote the core and far-field manifolds, respectively; this is equivalent to the notation used in previous works regarding localised radial patterns \cite{scheel2003radially,lloyd2009localized,mccalla2013spots,mcquighan2014oscillons}. Once both of these manifolds have been constructed, we evaluate their parametrisations at some coincidental point, which we can choose to be at $r=r_{0}$. Any function that lies on the intersection of the two manifolds is then, by definition, a localised solution.

Due to the explicit $r$-dependence of the system \eqref{Scheel:Syst}, the core manifold $\mathcal{W}^{cu}_{-}(\mu)$ must be constructed explicitly; this is usually done by introducing a variation-of-constants formula and applying a fixed-point argument. In contrast, because the radial coordinate $r$ is taken to be very large in the far-field region, the construction of the far-field manifold $\mathcal{W}^{s}_{+}(\mu)$ can be seen as an extension of standard local centre-manifold theory. In order to parametrise the far-field manifold $\mathcal{W}^{s}_{+}(\mu)$, it is helpful to reduce the problem onto the coordinates of a local centre manifold. In \cite{scheel2003radially}, it was proven that the extended vector field of \eqref{Scheel:Syst}, where $r$ is sufficiently large, can be reduced to the centre eigenspace; then, the reduced equations have a leading order expansion that is partially determined by the type of bifurcation occurring as $r\to\infty$. The two respective bifurcation types, homogeneous and Turing, are examined in \cite{scheel2003radially}; their reduced equations are found in \cite[Lemmas 3.9 \& 3.10]{scheel2003radially}. We now present one example of each bifurcation type, starting with a summary of the results of McQuighan and Sandstede \cite{mcquighan2014oscillons} regarding a homogeneous bifurcation.\par

\underline{\textit{Example 1: Homogeneous Bifurcation}} In \cite{mcquighan2014oscillons}, McQuighan and Sandstede investigate the existence of stationary localised radial patterns in the complex Ginzburg-Landau equation with $2:1$ forcing. By choosing the parameters of the system such that a homogeneous bifurcation occurs, the problem can be reduced to the real reaction-diffusion system,
\begin{align}
    0 = \Delta_{r}U - \mathbf{C}_{1}U - \varepsilon^{2}\mathbf{C}_{2}U - |U|^{2} \mathbf{C}_{3} U,\nonumber
\end{align}
where $U\in\mathbb{R}^{2}$, $\varepsilon$ is the bifurcation parameter, and $C_{i}$ is a square matrix that is determined by the parameters of the original system, for each $i=1,2,3$. This can also be written as the following first-order ODE system,
\begin{align}
    \mathbf{u}_{r} &= \mathcal{A}(r)\mathbf{u} + \mathcal{F}(\mathbf{u};\varepsilon),\qquad \mathcal{A}(r) = \begin{pmatrix}
    0 & 1_{N} \\ \mathbf{C}_{1} & -\frac{1}{r}\,1_{N}
    \end{pmatrix}, \qquad \mathcal{F}(\mathbf{u};\varepsilon) = \begin{pmatrix}
    0 \\ \varepsilon^{2}\mathbf{C}_{2}U + |U|^{2}\mathbf{C}_{3}U
    \end{pmatrix}, \label{McQuighan:Syst}
\end{align}
where $\mathbf{u}=(U,\partial_{r}U)^{\intercal}\in\mathbb{R}^{4}$. The linear solutions of \eqref{McQuighan:Syst} are associated with the eigenvectors of $\mathbf{C}_{1}$, which has eigenvalues $\lambda=0$ (double), $-m$, and $m$ for some $m\in\mathbb{R}$. For $\varepsilon>0$ the zero eigenvalues split onto the real line, as seen in Figure \ref{fig:Pulses} a), and so we have a homogeneous bifurcation as stated. By determining four linearly independent solutions of the linearisation of \eqref{McQuighan:Syst}, the core manifold $\mathcal{W}^{cu}_{-}(\varepsilon)$ was constructed using the standard approach of a fixed-point argument for a variation-of-constants formula; see \cite[Lemma 3.2]{mcquighan2014oscillons}. In order to construct the far-field manifold $\mathcal{W}^{s}_{+}(\varepsilon)$ and reduce the equations to some local centre manifold in the far-field region, McQuighan and Sandstede utilised the theory of strong stable foliations. This is presented in great detail in Section 4 of \cite{mcquighan2014oscillons}, and so we just note the following: there exists coordinates $A,B\in\mathbb{R}$ on some autonomous local centre manifold $\mathcal{W}^{c}_{+}(\varepsilon)$ such that \eqref{McQuighan:Syst} reduces to 
\begin{align}
    A_{r} &= B + \mathcal{R}_{A}(A,B,\sigma,\varepsilon),\nonumber\\
    B_{r} &= -\sigma A + \varepsilon^{2} A + c_{3}^{0}A^{3} + \mathcal{R}_{B}(A,B,\sigma,\varepsilon),\label{McQuighan:Normal}\\
    \sigma_{r} &= - \sigma^{2},\nonumber 
\end{align}
where $c_{3}^{0}$ is a fixed constant and $\mathcal{R}_{A}, \mathcal{R}_{B}$ contain higher-order terms. Here, $\sigma \propto \frac{1}{r}$ is taken to be an unconstrained variable such that the system becomes autonomous; in order to recover the $r$-dependent system \eqref{McQuighan:Syst}, we just define $\sigma(r_{0}) = \frac{1}{r_{0}}$ when evaluating the far-field parametrisation.

To find exponentially decaying solutions to \eqref{McQuighan:Normal}, McQuighan and Sandstede introduced geometric blow-up coordinates to track solutions backwards in $r$ to $r=r_{0}$. They determined that exponentially decaying solutions to \eqref{McQuighan:Normal} were associated with exponentially decaying solutions to the cubic non-autonomous Ginzburg-Landau equation,
\begin{align}
    \left(\partial_{ss} + \frac{1}{s}\partial_{s}\right)q(s)  = q(s) + c_{3}^{0}q(s)^{3}. \label{McQuighan:Ginzburg}
\end{align}
For $c_{3}^{0}<0$, there exists a unique exponentially decaying solution to \eqref{McQuighan:Ginzburg}; tracking back through the blow-up coordinate charts, McQuighan and Sandstede proved the existence of localised spikes in \eqref{McQuighan:Syst} such that the amplitude of the spike is $\textnormal{O}(\varepsilon)$.\par

\underline{\textit{Example 2: Turing Bifurcation}}
For the second example, we focus on the work of \cite{lloyd2009localized,mccalla2013spots} on the quadratic cubic Swift-Hohenberg equation,
\begin{align}
    0 = -\left( 1 + \Delta_{r}\right)^{2}u - \varepsilon^{2} u + \nu u^{2} - u^{3},\nonumber
\end{align}
where $u\in\mathbb{R}$, $\varepsilon$ is the bifurcation parameter, and $\nu$ is a fixed parameter. As with the previous examples, we express this equation as a first-order ODE system,
\begin{align}
    \begin{pmatrix}
    u_{1} \\ u_{2} \\ u_{3} \\ u_{4}
    \end{pmatrix}_{r} = \begin{pmatrix}
     0 & 0 & 1 & 0 \\
     0 & 0 & 0 & 1 \\
     -1 & 1 & -\frac{1}{r} & 0 \\
     0 & -1 & 0 & -\frac{1}{r} 
    \end{pmatrix}\begin{pmatrix}
    u_{1} \\ u_{2} \\ u_{3} \\ u_{4}
    \end{pmatrix} + \begin{pmatrix}
    0 \\ 0 \\ 0 \\ -\varepsilon^{2} u_{1} + \nu u_{1}^{2} - u_{1}^{3}
    \end{pmatrix}.\label{Swift:Syst}
\end{align}
Examining the spatial eigenvalues at the bifurcation point as $r\to\infty$, we observe a pair of purely imaginary eigenvalues at $\lambda = \pm\textnormal{i}$; these eigenvalues split off of the imaginary axis for $\varepsilon>0$, and so we have a Turing bifurcation as illustrated in Figure \ref{fig:Pulses} b). The core manifold $\mathcal{W}^{cu}_{-}(\mu)$ is constructed in the usual way; see \cite[Lemma 1]{lloyd2009localized}. The system \eqref{Swift:Syst} is already reduced to the centre eigenspace, and so no extra local centre manifold reduction is required. Introducing some complex coordinates $A,B\in\mathbb{C}$, \eqref{Swift:Syst} can be written as 
\begin{align}
    A_{r} &= - \frac{\sigma}{2}A + B + \mathcal{R}_{A}(A,B,\sigma,\varepsilon),\nonumber\\
    B_{r} &= -\frac{\sigma}{2} A + c_{0}\varepsilon^{2} A + c_{3}|A|^{2}A + \mathcal{R}_{B}(A,B,\sigma,\varepsilon),\label{Swift:Normal}\\
    \sigma_{r} &= - \sigma^{2},\nonumber 
\end{align}
where $c_{0},c_{3}$ are fixed constants and $\mathcal{R}_{A}, \mathcal{R}_{B}$ contain higher-order terms. Again, $\sigma \propto \frac{1}{r}$ is taken to be an unconstrained variable such that the system becomes autonomous; in order to recover the $r$-dependent system \eqref{Swift:Syst}, we just define $\sigma(r_{0}) = \frac{1}{r_{0}}$ when evaluating the far-field parametrisation.

To find exponentially decaying solutions to \eqref{Swift:Normal}, geometric blow-up coordinates are introduced to allow for tracking solutions backwards in $r$ to $r=r_{0}$. Exponentially decaying solutions to \eqref{Swift:Normal} were found to be associated with exponentially decaying solutions to the cubic non-autonomous Ginzburg-Landau equation,
\begin{align}
    \left(\partial_{s} + \frac{1}{2s}\right)^2 q(s)  = c_{0}q(s) + c_{3}q(s)^{3}. \label{Swift:Ginzburg}
\end{align}
For $c_{0}>0, c_{3}<0$, there exists a exponentially decaying solution to \eqref{Swift:Ginzburg}; there also exists a sufficiently small solution that decays exponentially independent of the sign of $c_{3}$. Tracking back through the blow-up coordinate charts, Lloyd and Sandstede \cite{lloyd2009localized} proved the existence of an elevated localised spot in \eqref{Swift:Syst}, as well as localised rings. Following this, McCalla and Sandstede \cite{mccalla2013spots} also proved the existence of a localised depressed spot, termed spot B. We will use these examples as a guide for our analysis of localised radial vegetation patterns.
\begin{figure}[t!]
    \centering
    \includegraphics[width=\textwidth]{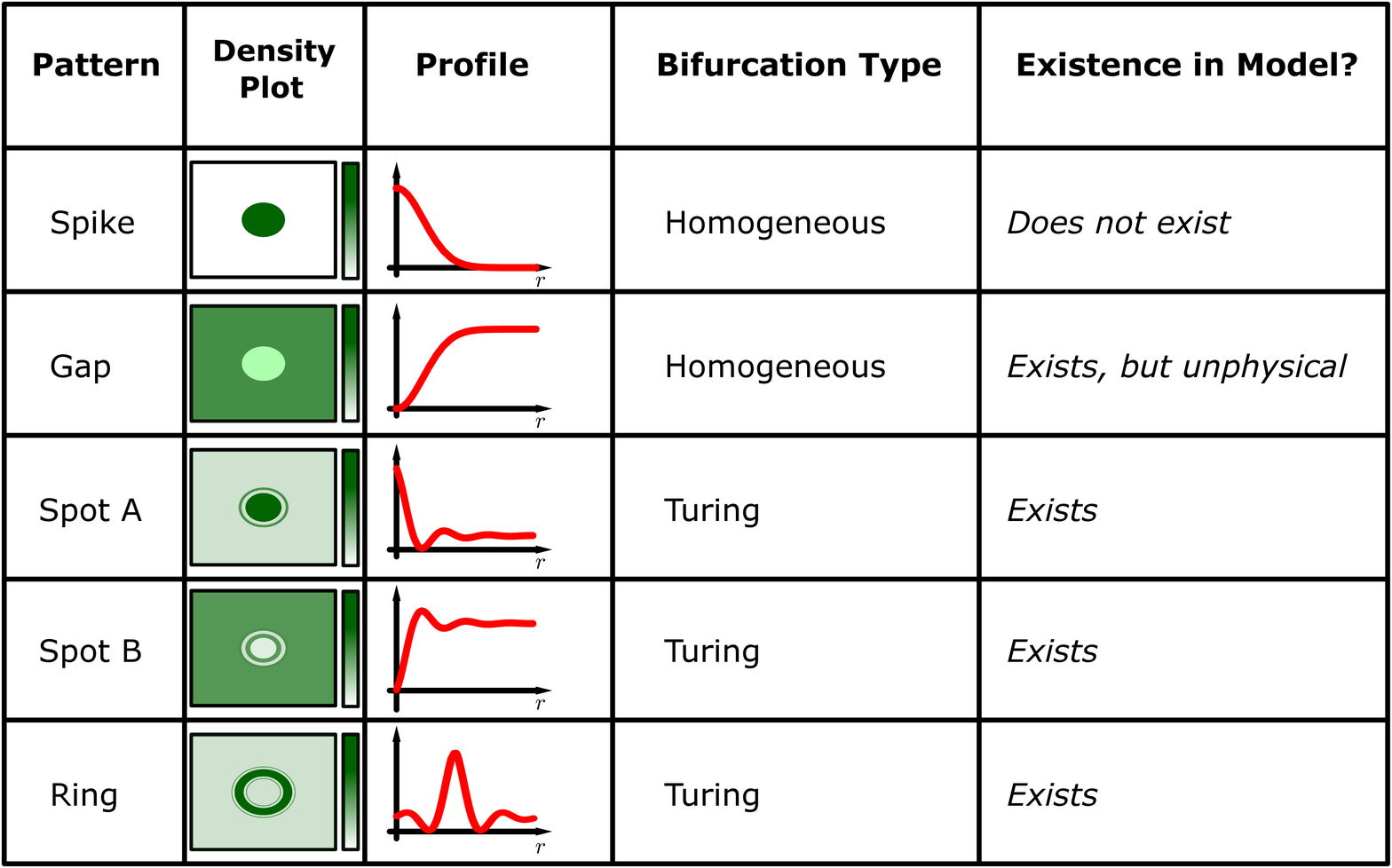}
    \caption{A summary of our results regarding the existence of different localised radial patterns. For the five classes of patterns, we present their relative vegetation density plot and radial profile, as well as their bifurcation type and whether or not they exist in the vegetation model.}
    \label{fig:results}
\end{figure}
\paragraph{Results} We summarise our results in Figure \ref{fig:results}. We look for two types of localised solutions bifurcating from the homogeneous instability: elevated spikes and depressed spikes, which we now refer to as `spikes' and `gaps', respectively. We find that no localised spikes bifurcate from the homogeneous instability of the von Hardenberg model, while localised gaps bifurcate both sub- and super-critically. However, for standard parameter values, localised gaps are seen to have a negative profile; hence, these solutions are unphysical. 

Close to a Turing instability, we prove that three classes of localised patterns exist and bifurcate sub-critically from a uniform vegetated state. These solutions include spots A and B, with respective maximum and minimum values at the core $r=0$, as well as elevated and depressed rings, whose maximum values are away from the core. We expect each of the localised patterns bifurcating from the Turing instability to remain physical for a reasonable choice of parameter values.

\paragraph{Overview} The rest of the paper is structured as follows: we begin by introducing the von Hardenberg model \cite{vonHardenberg2001} for stationary radial patterns, identifying uniform states and their associated bifurcation diagram. We then introduce the weakly nonlinear approximation for the von Hardenberg model derived in \cite{Dawes2016}; again, we compute uniform states and the bifurcation diagram of the system. We identify three parameter regions of the precipitation $p$ where localised states may emerge; these bifurcations are categorised by their spatial eigenvalues, and so we label each as either `homogeneous' or `Turing' bifurcations, appropriately.

Having categorised the bifurcations of the system, we relate each case to the relevant prototypical examples (\eqref{McQuighan:Syst} or \eqref{Swift:Syst}) and apply weakly nonlinear theory to compute the relevant coefficients required in the analysis. Where possible, we quote results from the prototypical examples; if the results are not identical we highlight the extra work necessary to complete the proof, although the technical details are beyond the scope of this paper.

Following this, we arrive at the following results for the reduced model: we prove the existence of three classes of localised radial patterns bifurcating from the uniform vegetated state at some Turing point $p=p_{1}$, where we recall that $p$ denotes precipitation. These localised patterns have spatially-oscillatory tails with an associated wave number $k_{1}>0$, and are directly associated with the spot A, spot B and ring solutions found in the Swift-Hohenberg equation \cite{lloyd2009localized,mccalla2013spots}. We also find evidence of the existence of localised gaps, characterised by depressions at their core and monotonic tails, bifurcating from the critical precipitation $p=p_{c}$ where the bare state destabilises. We expect there to be one solution bifurcating sub-critically from the bare state and another bifurcating super-critically from the uniform vegetated state, though these solutions may be unphysical.

We conclude by finding numerical solutions to the full von Hardenberg model corresponding to the localised spot A and localised gap solutions found in our analysis. We explore the effect of changing the shading parameter $\rho$ and plot the various bifurcation curves of these localised radial structures for distinct parameter choices. 

\section{The Model}
\label{s:model}
We focus our analysis on the model introduced by von Hardenberg et al.  \cite{vonHardenberg2001} to explain vegetation dynamics via a two-component reaction-diffusion system. This can be seen as an adaptation of the model presented in \cite{Klausmeier1999}; we note that there have been numerous studies of the von Hardenberg model (see \cite{Meron2004} for an in-depth description of model terms and behaviour, and \cite{Gowda2014} for an analysis of two-dimensional pattern selection), including the work of Dawes and Williams \cite{Dawes2016} in which localised one-dimensional patterns are observed. Hence, this model is a natural first choice for the study of localised radial patterns, since we can then directly compare our results with the one-dimensional case. We consider two non-negative dynamical variables $n(\mathbf{x},t)$ and $w(\mathbf{x},t)$, denoting the plant biomass density and ground water density, respectively; then, the dimensionless model takes the form,
\begin{align}
    n_{t} &=\left(\frac{\gamma w}{1+ \sigma w}\right)n - n^{2} - \nu n + \Delta n, \label{vonHard:1}\\
    w_{t} &= p - (1-\rho n)w - w^{2} n + \delta\Delta\left(w-\beta n\right) - c\left(w-\alpha n\right)_{x}. \nonumber
\end{align}
Here, $\Delta$ is the planar Laplacian, $p$ represents precipitation, and $\delta$ denotes the relative diffusivity of the water with respect to the vegetation, which is taken to be a large, positive parameter. In addition, $\gamma$, $\sigma$, $\rho$, $\beta$, and $\nu$ are taken to be $\textnormal{O}(1)$, positive parameters. We briefly summarise the ecological role of each term in \eqref{vonHard:1}. The term $(\frac{\gamma w}{1 + \sigma w})n$ models plant growth which has a linear rate for low levels of water. The terms $-\nu n$ and $-n^2$ account for the plant mortality and saturation due to limited nutrients, respectively, while the term $\Delta n$ models the spread of vegetation. The term $p$ denotes the increase of water due to precipitation, while the $-(1-\rho n)w$ term represents a loss of water due to evaporation. Hence, we require $n\leq \frac{1}{\rho}$ in order for the model to remain physical, where the parameter $\rho>0$ represents the reduction in evaporation caused by the presence of vegetation, mostly through shading effects. The term $-w^2 n$ models the local uptake of water, mostly determined by transpiration. The $\Delta(w - \beta n)$ term represents the transport of water in the soil as modelled by Darcy's law, where the $\beta n$ term accounts for the suction of water by the roots of the plant. The $-c(w-\alpha n)_{x}$ term models the advection of soil water down a slope in the $x$-direction, where $c$ is the runoff velocity. We will restrict our model to level ground, setting $c=0$, in order to study localised radial patterns. It would be interesting to reintroduce a slope, where $c> 0$, and investigate how localised radial patterns are affected by this imposed symmetry-breaking; this is left for future study.

\begin{figure}[t!]
    \centering
    \includegraphics[width=\textwidth]{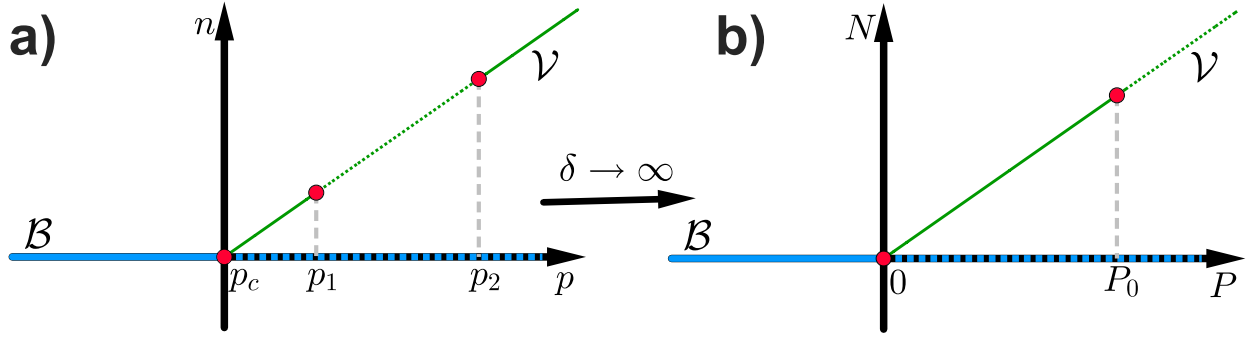}
    \caption{Diagram of the bifurcation structure for a) the full model \eqref{model:wn}, where the vegetation density $n$ is plotted against the precipitation $p$, and b) the weakly nonlinear reduction \eqref{model:NW}, where the scaled vegetation density $N$ is plotted against the scaled precipitation $P$. Stable solutions are indicated by solid lines and unstable solutions by dashed lines; points at which there is a change in stability are highlighted by red circles. In both equations, the uniform vegetated state $\mathcal{V}$ (green) bifurcates from the bare state $\mathcal{B}$ (blue) and undergoes at least one change of stability.}
    \label{fig:UniformStability}
\end{figure}

Setting $c=0$ and considering stationary axisymmetric solutions, \eqref{vonHard:1} becomes \begin{align}
    0 &=\left(\frac{\gamma w}{1+ \sigma w}\right)n - n^{2} - \nu n + \Delta_{r} n, \qquad \qquad 
    0 = -(w-p) + \left(\rho - w\right) w n + \delta\Delta_{r}\left(w-\beta n\right), \label{model:wn}
\end{align}
where $r$ is the planar radial coordinate and $\Delta_{r}:=\partial_{rr} + \frac{1}{r}\partial_{r}$ is the planar axisymmetric Laplacian. We begin by discussing the uniform states associated with \eqref{model:wn}.
\paragraph{Uniform States}
Any spatially uniform solution of \eqref{model:wn}, which we call $(n_{*},w_{*})$, must solve the following algebraic system,
\begin{align}
    0 &=\left[\left(\frac{\gamma w_{*}}{1+ \sigma w_{*}}\right) - n_{*} - \nu\right] n_{*}, \qquad \qquad 
    0 = -(w_{*}-p) + \left(\rho - w_{*}\right) w_{*} n_{*}. \nonumber
\end{align}
Hence, we see that there are two types of spatially uniform equilibria: the `bare' state $(n_{*}, w_{*}) = (0,p)$, and `vegetated' states where $n_{*} = \left(\frac{\gamma w_{*}}{1+ \sigma w_{*}}\right)-\nu$ and $w_{*}$ satisfies
\begin{align}
    -\left[\gamma - \sigma \nu\right]w_{*}^{3} + \left[\rho(\gamma - \sigma \nu) + \nu - \sigma\right]w_{*}^{2}  + \left[p \sigma - 1 - \rho\nu\right]w_{*} + p &= 0. \label{w*:eqn}
\end{align}
The bare state is known to undergo a homogeneous instability at $p=p_{c}:=\frac{\nu}{\gamma - \nu\sigma}$ (as discussed in \cite{Meron2004,Gowda2014,Dawes2016}) and so we require that $\gamma>\nu\sigma$ in order for $p_{c}$ to be physically relevant. We note that this condition also guarantees that \eqref{w*:eqn} has at least one positive solution $w_{*}$, by Descartes' rule of signs; for our choice of parameters, as seen in \cite{Gowda2014,Dawes2016}, there is only one uniform vegetated state that bifurcates from the Bare state at $p=p_{c}$. Turing bifurcations have been found at some values $p=p_{1}, p_{2}$, with associated critical wave numbers $k=k_{1}, k_{2}$, where periodic patterns emerge from the vegetated state; see \cite{Gowda2014}. For a bifurcation diagram of the uniform solutions to \eqref{model:wn}, see Figure \ref{fig:UniformStability} a). 

\paragraph{Weakly Nonlinear Reduction}
In order to make \eqref{model:wn} more analytically tractable, we now simplify the model in the case where $\delta\gg1$; to do this, we first introduce $\varepsilon:=\delta^{-\frac{1}{2}}$ such that we can write \eqref{model:wn} as
\begin{align}
    0 &=\left(\frac{\gamma w}{1+ \sigma w}\right)n - n^{2} - \nu n + \Delta_{r} n, \qquad \qquad 
    0 = -(w-p) + \left(\rho - w\right) w n + \frac{1}{\varepsilon^{2}}\Delta_{r}\left(w-\beta n\right). \label{model:eps}
\end{align}
This is now reminiscent of problems studied in \cite{Jaibi2020,Bastiaansen2019Stable,Carter2018}, in which existence of solutions can be proven via geometric singular perturbation methods. Similar techniques should be applicable here, although the $r$-dependent terms in \eqref{model:eps} may cause some extra complications. Instead, we will perform a weakly nonlinear reduction seen in \cite{Dawes2016} in order to find a simpler radial reaction-diffusion system that we can analyse. Following this, we will then numerically investigate the full system \eqref{model:wn} for physically relevant parameter values, using the intuition gained from the reduced model.

In \cite{Dawes2016}, Dawes and Williams presented a simplified model of \eqref{model:wn} for large values of $\delta$; for $\varepsilon=\delta^{-\frac{1}{2}}$, the authors found that
\begin{align}
    p_{1}= p_{c} + \textnormal{O}\left(\varepsilon^{2}\right),\qquad  w_{*}(p_{1})=w_{*}(p_{c})+\textnormal{O}\left(\varepsilon^{2}\right), \qquad n_{*}(p_{1})=\textnormal{O}\left(\varepsilon^{2}\right), \qquad k_{1}=\textnormal{O}\left(\varepsilon\right),\nonumber
\end{align}
where $(n_{*},w_{*})$ denotes the uniform vegetated state, $p_{1}$ is the minimum Turing point at which the vegetated state destabilises, and $k_{1}$ is the critical wave number associated with $p_{1}$. Then, by defining $p=p_{c} + \varepsilon^{2} P$, $n=\varepsilon^{2} N$, $w=p+\varepsilon^{2}W$ and $\widetilde{r}=\varepsilon r$, \eqref{model:eps} can be rewritten, to leading order, as 
\begin{align}
    0 &=\Delta_{\widetilde{r}} N + b_{0}P N + \left(b_{0}W - N\right)N, \qquad \qquad 
    0 = \Delta_{\widetilde{r}}\left(W-\beta N\right) - W + a_{0}N \label{model:NW}
\end{align}
as $\varepsilon\to0$, where $b_{0}= \frac{\gamma }{\left(1 + \sigma p_{c}\right)^{2}}$, $p_{c}= \frac{\nu}{\gamma - \sigma\nu}$, and $a_{0}=\left(\rho - p_{c}\right)p_{c}$. For simplicity, we perform the transformation $\left(N,W,P\right)\mapsto\left(N,\beta W, b_{0}^{-1}P\right)$ and drop the tilde notation for $r$, such that we can write \eqref{model:NW} as
\begin{align}
    \Delta_{r}\begin{pmatrix}N \\ W\end{pmatrix} &= \begin{pmatrix} -P N - \left(b W - N\right)N \\ W - a N -P N - \left(b W - N\right)N\end{pmatrix},\label{model:N,W}
\end{align}
where $a=\frac{a_{0}}{\beta}$ and $b=b_{0}\beta$. Again, we see that there are two uniform states: the `bare' state $\mathcal{B}:(N,W)=(0,0)$, and a uniform vegetated state $\mathcal{V}:(N,W)=(N_{0},W_{0})$, where
\begin{align}
     N_{0} = \frac{P}{1 - a b}, \qquad \qquad W_{0} = \frac{a P}{1 - a b}.
     \label{Unif:V}
 \end{align} 
 Then, we require $\textnormal{sgn}(P)=\textnormal{sgn}(1-a b)$ such that $N_{0}>0$. We note that $1 > a b$ for the choice of parameters seen in \cite{Gowda2014,Dawes2016}; see \eqref{param:Gowda}, \eqref{param:Dawes}. Hence, we see that the vegetated state is bifurcating from the bare state super-critically; the bifurcation diagram for the uniform states of \eqref{model:N,W} is presented in Figure \ref{fig:UniformStability} b).
\section{Localised Radial Patterns in the Reduced Model}
We now wish to identify possible bifurcations of the reduced system \eqref{model:N,W}, where localised solutions may emerge as perturbations from a uniform state. To do this, we will linearise about each uniform state in turn, investigating the spatial eigenvalues in the limit as $r\to\infty$.
\subsection{Bifurcation Analysis and Canonical Forms}
\label{s:Bif}
\begin{figure}[t!]
    \centering
    \includegraphics[width=\textwidth]{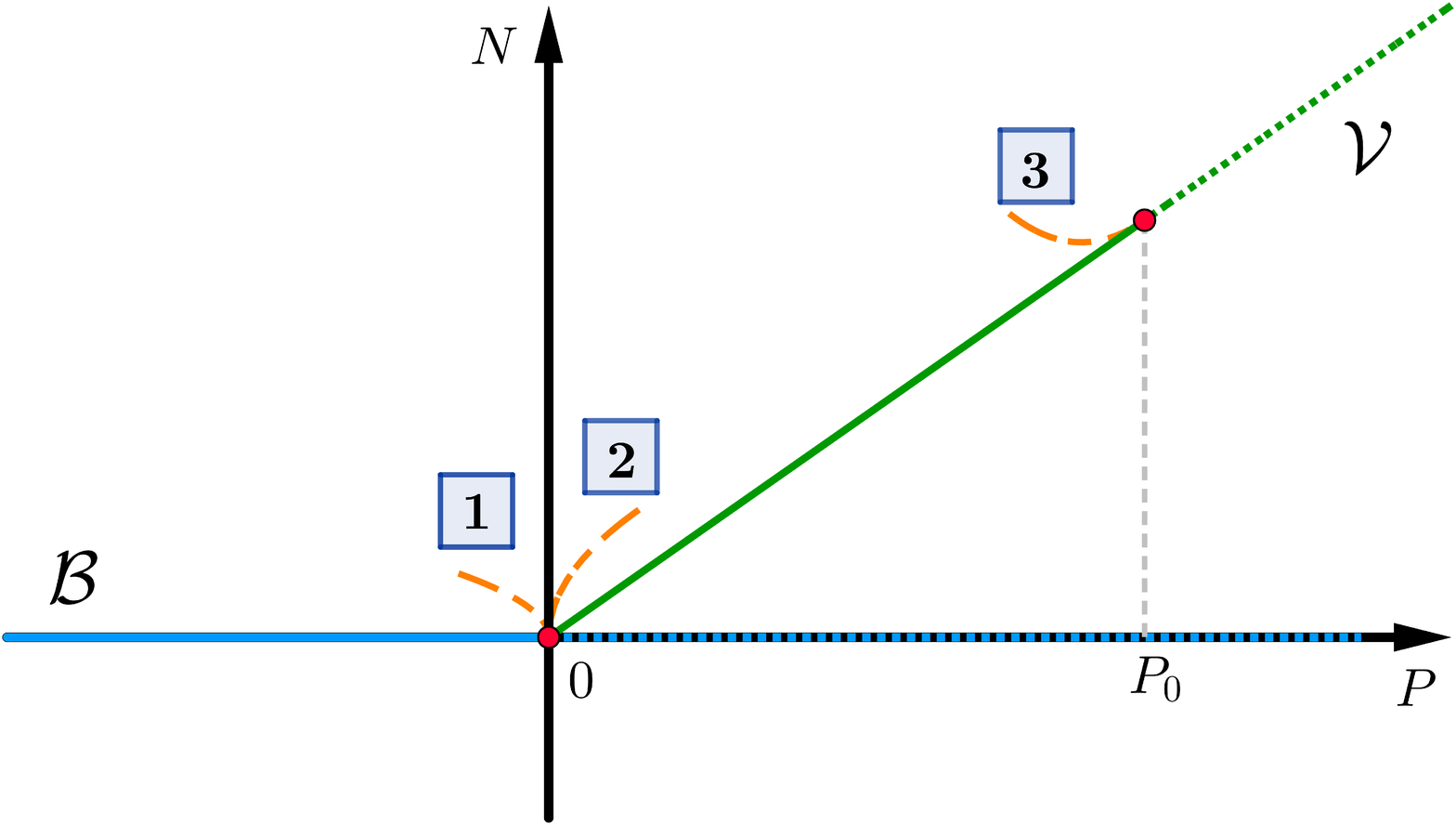}
    \caption{In the weakly reduced model \eqref{model:NW}, we expect there to be three possible cases in which localised radial solutions can emerge. In case 1, solutions bifurcate from the bare state $\mathcal{B}$ (blue) sub-critically at $P=0$. In case 2, solutions bifurcate from the uniform vegetated state $\mathcal{V}$ (green) super-critically at $P=0$. In case 3, solutions bifurcate from the uniform vegetated state $\mathcal{V}$ (green) sub-critically at $P=P_{0}$.}
    \label{fig:BifurcationPoints}
\end{figure}

We recall that there are two uniform states associated with \eqref{model:N,W}: the bare state $\mathcal{B}:(N,W)=(0,0)$, and the vegetated state $\mathcal{V}:(N,W)=(N_{0},W_{0})$, where $N_{0},W_{0}$ are defined in \eqref{Unif:V}. We begin by looking for any localised radial solutions bifurcating from the bare state $\mathcal{B}$.
\paragraph{Case 1}
 We linearise about the bare state $\mathcal{B}$; then, we can express \eqref{model:N,W} as the first-order system 
\begin{align}
 \mathbf{U}_{r} &= \mathbf{L}_{\mathcal{B}}\left(r\right)\mathbf{U} + \mathbf{F}_{\mathcal{B}}\left(\mathbf{U}\right),\label{model:U;B}
\end{align}
 for $\mathbf{U}=\left(N,W,\partial_{r}N,\partial_{r}W\right)\in\mathbb{R}^{4}$, where
 \begin{align}
     \mathbf{L}_{\mathcal{B}}\left(r\right) = \begin{pmatrix}0 & 0 & 1 & 0\\ 0 & 0 & 0 & 1\\ -P & 0 & -\frac{1}{r} & 0\\ -P-a & 1 & 0 & -\frac{1}{r} \end{pmatrix},\qquad \qquad \mathbf{F}_{\mathcal{B}}\left(\mathbf{U}\right) = \begin{pmatrix}0 \\ 0 \\ - \left(b \mathbf{U}_{2} - \mathbf{U}_{1}\right)\mathbf{U}_{1} \\ - \left(b \mathbf{U}_{2} - \mathbf{U}_{1}\right)\mathbf{U}_{1} \end{pmatrix}.\nonumber
 \end{align}
 We investigate the eigenvalues of the linear operator $\mathbf{L}_{\mathcal{B}}(r)$ in the limit as $r\to\infty$; any eigenvalue $\lambda\in\mathbb{C}$ of $\lim_{r\to\infty}\mathbf{L}_{\mathcal{B}}(r)$ must be a root of the characteristic polynomial,
\begin{align}
    p_{\mathcal{B}}(\lambda) &= \lambda^{4} - \left(1-P\right)\lambda^{2} - P = \left(\lambda^{2}+P\right)\left(\lambda^{2}-1\right).\nonumber
\end{align}
We see that solutions of $p_{\mathcal{B}}(\lambda)=0$ take the form $\lambda=1,-1$, and $\pm\sqrt{-P}$. Hence, for $P>0$ we have purely imaginary eigenvalues $\pm \textnormal{i}\sqrt{P}$; these eigenvalues pass through the origin at $P=0$ and, for $P<0$, lie on the real axis at $\pm\sqrt{|P|}$. Thus, we define $P=-\varepsilon^{2}$ such that the bifurcation point is at $\varepsilon=0$, and we expect localised solutions to emerge in the region $0<\varepsilon\ll1$. Taking $\varepsilon$ small, we now write \eqref{model:U;B} as
\begin{align}
 \mathbf{U}_{r} &= \mathbf{L}_{1}\left(r\right)\mathbf{U} + \mathbf{F}_{1}\left(\mathbf{U};\varepsilon^{2}\right),\label{model:1}
\end{align}
 where
 \begin{align}
     \mathbf{L}_{1}\left(r\right) = \begin{pmatrix}0 & 0 & 1 & 0\\ 0 & 0 & 0 & 1\\ 0 & 0 & -\frac{1}{r} & 0\\ -a & 1 & 0 & -\frac{1}{r} \end{pmatrix},\qquad \qquad \mathbf{F}_{1}\left(\mathbf{U};\varepsilon^{2}\right) = \begin{pmatrix}0 \\ 0 \\ \varepsilon^{2}\mathbf{U}_{1}- \left(b \mathbf{U}_{2} - \mathbf{U}_{1}\right)\mathbf{U}_{1} \\ \varepsilon^{2}\mathbf{U}_{1}- \left(b \mathbf{U}_{2} - \mathbf{U}_{1}\right)\mathbf{U}_{1} \end{pmatrix}.\nonumber
 \end{align}
 This is the only bifurcation for the bare state $\mathcal{B}$, and so we now turn our attention to localised solutions bifurcating off of the vegetated state $\mathcal{V}$.
\paragraph{Case 2}
We linearise about the vegetated state $\mathcal{V}$; that is, we define $N=N_{0} + N_{1}$, $W=W_{0} + W_{1}$ where $N_{0}$ and $W_{0}$ are defined in \eqref{Unif:V} and $N_{1}, W_{1}$ are taken to be small. Then, we can write \eqref{model:N,W} as the first-order system,
\begin{align}
 \mathbf{V}_{r} &= \mathbf{L}_{\mathcal{V}}\left(r\right)\mathbf{V} + \mathbf{F}_{\mathcal{V}}\left(\mathbf{V}\right),\label{model:U;V}
\end{align}
 for $\mathbf{V}=\left(N_{1},W_{1},\partial_{r}N_{1},\partial_{r}W_{1}\right)\in\mathbb{R}^{4}$, where
 \begin{align}
     \mathbf{L}_{\mathcal{V}}\left(r\right) = \begin{pmatrix}0 & 0 & 1 & 0\\ 0 & 0 & 0 & 1\\ \frac{P}{1- a b} & -\frac{b P}{1- a b} & -\frac{1}{r} & 0\\ \frac{P}{1- a b}-a & 1 - \frac{b P}{1- a b} & 0 & -\frac{1}{r} \end{pmatrix},\qquad \qquad \mathbf{F}_{\mathcal{V}}\left(\mathbf{V}\right) = \begin{pmatrix}0 \\ 0 \\ - \left(b \mathbf{V}_{2} - \mathbf{V}_{1}\right)\mathbf{V}_{1} \\ - \left(b \mathbf{V}_{2} - \mathbf{V}_{1}\right)\mathbf{V}_{1} \end{pmatrix}.\nonumber
 \end{align}
 We investigate the eigenvalues of the linear operator $\mathbf{L}_{\mathcal{V}}(r)$ in the limit as $r\to\infty$; an eigenvalue $\lambda\in\mathbb{C}$ of $\lim_{r\to\infty}\mathbf{L}_{\mathcal{V}}(r)$ must be a root of the characteristic polynomial,
\begin{align}
    p_{\mathcal{V}}(\lambda) &= \lambda^{4} - \left(1 + K P\right)\lambda^{2} + P,\label{p:V}
\end{align}
where we have defined $K:=\frac{1-b}{1-a b}$. For small values of $P$, we see that 
\begin{align}
    \lambda &= 1 + \left(\frac{K-1}{2}\right)P + \textnormal{O}\left(P^{2}\right), \quad \lambda = -1 - \left(\frac{K-1}{2}\right)P + \textnormal{O}\left(P^{2}\right), \quad \textnormal{or} \quad  \lambda = \pm\sqrt{P} + \textnormal{O}\left(P^{\frac{3}{2}}\right).\nonumber
\end{align}
Hence, taking $P= \varepsilon^{2}$, we again have a bifurcation at $\varepsilon=0$ such that we expect localised solutions to emerge in the region $0<\varepsilon\ll1$. We can now write the roots of \eqref{p:V} as 
\begin{align}
    \lambda &= 1 + \textnormal{O}(\varepsilon), \quad \lambda = -1 + \textnormal{O}(\varepsilon), \quad \lambda = \varepsilon + \textnormal{O}(\varepsilon^{3}), \quad\textnormal{or}\quad  \lambda = -\varepsilon + \textnormal{O}(\varepsilon^{3}),\nonumber
\end{align}
and express \eqref{model:U;V} as
\begin{align}
 \mathbf{V}_{r} &= \mathbf{L}_{1}\left(r\right)\mathbf{V} + \mathbf{F}_{2}\left(\mathbf{V};\varepsilon^{2}\right),\label{model:2}
\end{align}
 where $\mathbf{L}_{1}(r)$ is the linear operator defined in \eqref{model:1}, and 
 \begin{align}
     \mathbf{F}_{2}\left(\mathbf{V};\varepsilon^{2}\right) = \begin{pmatrix}0 \\ 0 \\ \varepsilon^{2}\left(\frac{\mathbf{V}_{1} - b \mathbf{V}_{2}}{1- a b } \right) - \left(b \mathbf{V}_{2} - \mathbf{V}_{1}\right)\mathbf{V}_{1} \\ \varepsilon^{2}\left(\frac{\mathbf{V}_{1} - b \mathbf{V}_{2}}{1- a b } \right) - \left(b \mathbf{V}_{2} - \mathbf{V}_{1}\right)\mathbf{V}_{1} \end{pmatrix}.\nonumber
 \end{align}
Thus, at the bifurcation point $P=0$, we have a homogeneous bifurcation close to which localised radial solutions may emerge, bifurcating from either uniform state. Notably, any solutions bifurcating from the bare state $\mathcal{B}$ will emerge in a sub-critical region ($P<0$), whereas any solutions bifurcating from the vegetated state $\mathcal{V}$ will emerge in a super-critical region ($P>0$). As previously discussed, we also expect the vegetated state $\mathcal{V}$ to undergo a Turing instability, and so we return to our characteristic polynomial $p_{\mathcal{V}}(\lambda)$.
\paragraph{Case 3}
For the vegetated state $\mathcal{V}$, we expect a Turing instability to be associated with the repeated roots of $p_{\mathcal{V}}(\lambda)$, as defined in \eqref{p:V}. Completing the square, we note that
 \begin{align}
    p_{\mathcal{V}}(\lambda) &=  \left(\lambda^{2} - \frac{1}{2}\left[1 + K P\right]\right)^{2} + P - \frac{1}{4}\left(1 + K P\right)^{2},\nonumber
\end{align}
and so $p_{\mathcal{V}}(\lambda)$ has repeated roots at some critical $P=P_{c}$, where $P_{c}$ satisfies
 \begin{align}
     0 &= 4 P_{c} - \left(1 + K P_{c}\right)^{2},\qquad \qquad 
     P_{c,\pm} = \frac{2-K\pm2\sqrt{1-K}}{K^{2}}.\nonumber
\end{align}
 Both choices of $P_{c, \pm}$ have an associated eigenvalue $\lambda_{\pm}$, where
 \begin{align}
    \lambda_{\pm}^{2} &= \frac{1}{2}\left[1 + K P_{c,\pm}\right] = \frac{1\pm\sqrt{1-K}}{K}.\nonumber
\end{align}
For the parameter choices seen in \cite{Gowda2014,Dawes2016}, we find that $K<0$; hence, $\lambda_{+}^{2}<0$ and $\lambda_{-}^{2}>0$, and so a Turing instability occurs at $P=P_{0}:=\frac{2-K + 2\sqrt{1-K}}{K^{2}}$. The critical eigenvalue is defined as $\lambda = \pm \textnormal{i}k$, where $k := \sqrt{\frac{1\pm\sqrt{1-K}}{-K}}\in\mathbb{R}$. We note that $P_{0} = k^{4}$ and take $P= k^{4} + \widetilde{P}$; then, we see that
\begin{align}
    \lambda^{2} &= \frac{1}{2}\left[1 - \left(\frac{2k^{2}+1}{k^{4}}\right) \left(k^{4} + \widetilde{P}\right)\right] \pm \sqrt{\frac{1}{4}\left[1 - \left(\frac{2k^{2}+1}{k^{4}}\right) \left(k^{4} + \widetilde{P}\right)\right]^{2} - \left(k^{4} + \widetilde{P}\right)},\nonumber\\
    &= -\left[k^{2} + \left(\frac{2k^{2}+1}{2k^{4}}\right)\widetilde{P}\right] \pm \sqrt{\left(\frac{2k^{2}+1}{2k^{4}}\right)^{2}\widetilde{P} + \left(\frac{k^{2}+1}{k^{2}}\right)} \sqrt{\widetilde{P}}, \nonumber\\
    &= -k^{2} \pm \sqrt{\left(\frac{k^{2}+1}{k^{2}}\right)}\sqrt{\widetilde{P}} - \left(\frac{2k^{2}+1}{2k^{4}}\right)\widetilde{P}+\textnormal{O}\left(\widetilde{P}^{\frac{3}{2}}\right).\nonumber
\end{align}
We recall that the Turing instability occurs as the spatial eigenvalues bifurcate off of the imaginary axis, as seen in Figure \ref{fig:Pulses} b). This would mean that $\lambda^{2}$ has non-zero imaginary part, and so $\widetilde{P}$ must be negative; that is, we expect the Turing instability to be a sub-critical bifurcation. We define $P:=P_{0} - \varepsilon^{2}$ such that $\varepsilon=0$ is the bifurcation point and we expect localised solutions to emerge for $0<\varepsilon\ll1$. Hence, we can now write the eigenvalues $\{\lambda_{i}\}_{i=1}^{4}$ as
\begin{align}
    \lambda_{1} &= \textnormal{i}k + \textnormal{O}\left(\varepsilon\right), \qquad \lambda_{2} = \textnormal{i}k + \textnormal{O}\left(\varepsilon\right),\qquad
    \lambda_{3} = -\textnormal{i}k + \textnormal{O}\left(\varepsilon\right), \qquad \lambda_{4} = -\textnormal{i}k + \textnormal{O}\left(\varepsilon\right).\nonumber
\end{align}
and express \eqref{model:U;V} as
\begin{align}
 \mathbf{V}_{r} &= \mathbf{L}_{3}\left(r\right)\mathbf{V} + \mathbf{F}_{3}\left(\mathbf{V};\varepsilon^{2}\right),\label{model:3}
\end{align}
 where
 \begin{align}
     \mathbf{L}_{3}\left(r\right) = \begin{pmatrix}0 & 0 & 1 & 0\\ 0 & 0 & 0 & 1\\ d_{k} - k^{2} & -\left[1 + k^{2} + d_{k}\right] & -\frac{1}{r} & 0\\ \frac{d_{k}^{2}}{1+k^{2}+d_{k}} & -d_{k} - k^{2} & 0 & -\frac{1}{r} \end{pmatrix},\qquad  \mathbf{F}_{3}\left(\mathbf{V};\varepsilon^{2}\right) = \begin{pmatrix}0 \\ 0 \\ \frac{d_{k}-k^{2}}{k^{4}}\left(b \mathbf{V}_{2} - \mathbf{V}_{1}\right)\varepsilon^{2} - \left(b \mathbf{V}_{2} - \mathbf{V}_{1}\right)\mathbf{V}_{1} \\ \frac{d_{k}-k^{2}}{k^{4}}\left(b \mathbf{V}_{2} - \mathbf{V}_{1}\right)\varepsilon^{2} - \left(b \mathbf{V}_{2} - \mathbf{V}_{1}\right)\mathbf{V}_{1} \end{pmatrix},\nonumber
 \end{align}
 and we have defined $d_{k}:=\frac{k^{4}}{1-ab} + k^{2} =  \left(\frac{2k^{2}+1}{b-1}\right) + k^{2}$. We introduce the following transformation
\begin{align}
     \mathbf{P}:=\begin{pmatrix}
     \mathcal{P} & \mathbf{0} \\ \mathbf{0} & \frac{1}{k}\mathcal{P}
     \end{pmatrix}, \qquad \qquad \mathcal{P}:=\begin{pmatrix}
     \frac{1}{(d_{k}-k^{2})} & 0\\
     \frac{d_{k}}{k^{2}(d_{k}-k^{2})} & -\frac{d_{k}+k^{2}+1}{k^{2}(d_{k}-k^{2})}\end{pmatrix},\label{Defn:P}
 \end{align}
 which has a well-defined inverse
 \begin{align}
     \mathbf{P}^{-1}:=\begin{pmatrix}
     \mathcal{P}^{-1} & \mathbf{0} \\ \mathbf{0} & k\mathcal{P}^{-1}
     \end{pmatrix}, \qquad \qquad \mathcal{P}^{-1}:=\begin{pmatrix}
     (d_{k}-k^{2}) & 0\\
     \frac{d_{k}(d_{k}-k^{2})}{d_{k} + k^{2} + 1} & -\frac{k^{2}(d_{k}-k^{2})}{d_{k} + k^{2} + 1}\end{pmatrix}.\nonumber
 \end{align}
 Then, defining $\widetilde{\mathbf{V}}:=\mathbf{P}\mathbf{V}$, $\widetilde{r}:=k r$ and $\widetilde{\varepsilon}:=k^{-2}\varepsilon$, we can finally write \eqref{model:3} as 
 \begin{align}
 \widetilde{\mathbf{V}}_{\widetilde{r}} &= \widetilde{\mathbf{L}}_{3}\left(\widetilde{r}\right)\widetilde{\mathbf{V}} + \widetilde{\mathbf{F}}_{3}\left(\widetilde{\mathbf{V}};\widetilde{\varepsilon}^{2}\right),\label{model:3;tilde}
\end{align}
 where
 \begin{align}
     \widetilde{\mathbf{L}}_{3}\left(\widetilde{r}\right) = \begin{pmatrix}0 & 0 & 1 & 0\\ 0 & 0 & 0 & 1\\ -1 & 1 & -\frac{1}{\widetilde{r}} & 0\\ 0 & -1 & 0 & -\frac{1}{\widetilde{r}} \end{pmatrix},\qquad  \widetilde{\mathbf{F}}_{3}\left(\widetilde{\mathbf{V}};\widetilde{\varepsilon}^{2}\right) = \begin{pmatrix}0 \\ 0 \\ \left( \widetilde{\mathbf{V}}_{1} -  \widetilde{\mathbf{V}}_{2}\right)\left(\widetilde{\varepsilon}^{2} -\widetilde{\mathbf{V}}_{1}\right) \\ -\omega\left( \widetilde{\mathbf{V}}_{1} -  \widetilde{\mathbf{V}}_{2}\right)\left(\widetilde{\varepsilon}^{2} -\widetilde{\mathbf{V}}_{1}\right) \end{pmatrix},\nonumber
 \end{align}
with the fixed parameter $\omega:=\frac{k^{2}+1}{k^{2}}$. Finally, we drop the tilde notation for simplicity. Thus, we have three regions of interest: the sub- and super-critical regions close to a homogeneous bifurcation at $P=0$, where localised solutions may bifurcate from the respective bare or vegetated states, and the sub-critical region close to a Turing bifurcation at $P=P_{0}$, where localised solutions may bifurcate from the vegetated state; see Figure \ref{fig:BifurcationPoints}. We now look to relate these cases to the prototypical equations seen in \cite{mcquighan2014oscillons,lloyd2009localized,mccalla2013spots} associated with each bifurcation type, and then discuss the existence of possible localised radial solutions in each case.
\subsection{Cases 1 and 2: Homogeneous Bifurcation}
\label{s:Cases1&2}
We recall that cases 1 and 2 both occur close to a  homogeneous bifurcation, and so we expect that we can relate both to the study of oscillons in the complex Ginzburg-Landau equation \cite{mcquighan2014oscillons}. To make this more clear, we express our equations \eqref{model:1}, \eqref{model:2} in a more general form; that is, for cases 1 and 2, we write our equations in the following way,
\begin{align}
 \begin{pmatrix}
 \mathbf{U} \\\mathbf{V}
 \end{pmatrix}_{r} &= \begin{pmatrix}
 \mathbf{0} & 1_{2} \\\mathbf{C} & -\frac{1}{r}1_{2}
 \end{pmatrix}\begin{pmatrix}
 \mathbf{U} \\\mathbf{V}
 \end{pmatrix} + \begin{pmatrix}
 \mathbf{0} \\\varepsilon^{2}\mathbf{D}\mathbf{U} + \mathbf{f}(\mathbf{U})
 \end{pmatrix}, \label{model:gen}
\end{align}
 for some $\mathbf{U},\mathbf{V}\in\mathbb{R}^{2}$, $\mathbf{C},\mathbf{D}\in\mathbb{R}^{2\times2}$, and $\mathbf{f}:\mathbb{R}^{2}\to\mathbb{R}^{2}$. In both cases, 
 \begin{align}
     \mathbf{C} &= \begin{pmatrix}
 0 & 0 \\ -a & 1
 \end{pmatrix}, \qquad \qquad \mathbf{f}:\begin{pmatrix}
 U_{1} \\ U_{2}
 \end{pmatrix}\mapsto \begin{pmatrix}
 -(b\,U_{2} - U_{1})U_{1} \\ -(b\,U_{2} - U_{1})U_{1}
 \end{pmatrix},\nonumber\\
 \intertext{however $\mathbf{U}$, $\mathbf{V}$ and $\mathbf{D}$ differ for each case. In particular, for case 1,}
      (1):\qquad\mathbf{U} &=\begin{pmatrix}
 N \\W
 \end{pmatrix}, \qquad \mathbf{V}=\begin{pmatrix}
 \partial_{r}N \\ \partial_{r}W
 \end{pmatrix}, \qquad \mathbf{D} = \begin{pmatrix}
 1 & 0 \\ 1 & 0
 \end{pmatrix},\nonumber\\
 \intertext{whereas in case 2, we have}
 (2):\qquad\mathbf{U} &=\begin{pmatrix}
 N_{1} \\W_{1}
 \end{pmatrix}, \qquad \mathbf{V}=\begin{pmatrix}
 \partial_{r}N_{1} \\ \partial_{r}W_{1}
 \end{pmatrix}, \qquad \mathbf{D} = \begin{pmatrix}
  \frac{1}{1-ab} & -\frac{b}{1-ab} \\ \frac{1}{1-ab} &  -\frac{b}{1-ab}
 \end{pmatrix}.\nonumber
 \end{align}
 By writing our equations in this form, we have expressed cases 1 and 2 in an analogous way to the CGL problem in \cite{mcquighan2014oscillons}. We note that the general system \eqref{model:gen} has a quadratic nonlinearity, which is in contrast to $\mathbb{Z}_{2}$ symmetry seen in \cite{mcquighan2014oscillons}. However, the core parametrisation and the foliations in the far-field are entirely dependent on the linear part of the system, and so the results in \cite{mcquighan2014oscillons} should follow identically if $\mathbf{C}$ has the same properties as the matrix $C_{1}$ in \cite{mcquighan2014oscillons}. Hence, we begin by analysing the associated linear problem of \eqref{model:gen}.
  
We note that, at the bifurcation point $\varepsilon=0$, the linear operator of \eqref{model:gen} in the limit as $r\to\infty$ has eigenvalues $\lambda^{c}=0$ (with double multiplicity), $\lambda^{u}=1$, and $\lambda^{s}=-1$; these eigenvalues have associated generalised eigenspaces
 \begin{align}
     E^{c}_{+}=\textnormal{span}\,\left\{\begin{pmatrix}\widetilde{U}_{0} \\ 0\end{pmatrix}, \begin{pmatrix}0 \\ \widetilde{U}_{0}\end{pmatrix}\right\},\qquad E^{u}_{+}=\textnormal{span}\,\left\{\begin{pmatrix}\widetilde{U}_{1}\\ \widetilde{U}_{1}\end{pmatrix}\right\},\qquad E^{s}_{+}=\textnormal{span}\,\left\{\begin{pmatrix}\widetilde{U}_{1} \\ -\widetilde{U}_{1}\end{pmatrix}\right\},\nonumber\end{align}
     where,
    \begin{align}
     \widetilde{U}_{0} = \begin{pmatrix} 1 \\ a\end{pmatrix}, \qquad \widetilde{U}_{1} = \begin{pmatrix} 0 \\ 1\end{pmatrix}. \label{Defn:U0}
 \end{align}
We note that this has the same structure as in \cite{mcquighan2014oscillons} (in the case when $m=1$) for different values of $\widetilde{U}_{0},\widetilde{U}_{1}$. Then, the linearisation of \eqref{model:gen} at $\varepsilon=0$,
\begin{align}
 \begin{pmatrix}
 \mathbf{U} \\\mathbf{V}
 \end{pmatrix}_{r} &= \begin{pmatrix}
 \mathbf{0} & 1_{2} \\\mathbf{C} & -\frac{1}{r}1_{2}
 \end{pmatrix}\begin{pmatrix}
 \mathbf{U} \\\mathbf{V}
 \end{pmatrix}, \label{model:1;lin}
\end{align}
for some $\mathbf{U},\mathbf{V}\in\mathbb{R}^{2}$, has four linearly independent solutions $\{\mathbf{W}_{i}(r)\}_{i=1}^{4}$, given by
\begin{align}
    \mathbf{W}_{1}(r) = \begin{pmatrix} \widetilde{U}_{0} \\ 0\end{pmatrix},\qquad \mathbf{W}_{2}(r) = \begin{pmatrix} \widetilde{U}_{0} \ln(r)\\ \widetilde{U}_{0} \frac{1}{r}\end{pmatrix},\qquad \mathbf{W}_{3}(r) = \begin{pmatrix} \widetilde{U}_{1} I_{0}(r) \\ \widetilde{U}_{1} I_{1}(r) \end{pmatrix},\qquad \mathbf{W}_{4}(r) = \begin{pmatrix} \widetilde{U}_{1} K_{0}(r) \\ -\widetilde{U}_{1} K_{1}(r) \end{pmatrix},\label{Solns:1;lin}
\end{align}
 where $I_{\nu}(r),K_{\nu}(r)$ are the $\nu^{th}$-order modified Bessel functions of the first and second kind, respectively; see \cite[Section 9.6]{abramowitz1972handbook}. We also note that $\textnormal{span}\,\{\mathbf{W}_{1}(r), \mathbf{W}_{2}(r)\} \equiv E^{c}_{+}$ for all $r$, since $c_{1}\mathbf{W}_{1}(r) + c_{2}\mathbf{W}_{2}(r) = \left(c_{1} + c_{2} \ln (r)\right)\mathbf{W}^{0}_{1} + \frac{c_{2}}{r}\mathbf{W}^{0}_{2}$, where $\mathbf{W}^{0}_{1}:=\left(\widetilde{U}_{0},\mathbf{0}\right)^{\intercal}$ and $\mathbf{W}^{0}_{2}:=\left(\mathbf{0},\widetilde{U}_{0}\right)^{\intercal}$. We define $P_{-}^{cu}(r)$ to be the projection on $\textnormal{span}\,\{\mathbf{W}_{1}(r), \mathbf{W}_{3}(r)\}$ along $\textnormal{span}\,\{\mathbf{W}_{2}(r), \mathbf{W}_{4}(r)\}$ for each fixed $r\geq0$.
 
 These solutions have exactly the same structure as \cite{mcquighan2014oscillons}; since the explicit values of $\widetilde{U}_{0}$ and $\widetilde{U}_{1}$ are not used in the following analysis, we can conclude the following results from \cite{mcquighan2014oscillons}.
 
 The first result is regarding the existence and parametrisation of the core manifold $\mathcal{W}_{-}^{cu}(\varepsilon)$, which contains all small-amplitude solutions of \eqref{model:gen} on the interval $r\in[0,r_{0}]$, for some large fixed $r_{0}>0$, that remain bounded as $r\to0$,
 \begin{Lemma}
  Fix $r_{0}>0$; then, there are constants $\delta_{0},\delta_{1},\varepsilon_{0}>0$ such that the set $\mathcal{W}^{cu}_{-}(\varepsilon)$ of solutions $\mathbf{U}(r)$ of \eqref{model:gen} for which $\sup_{0\leq r\leq r_{0}}|\mathbf{U}|<\delta_{0}$ is a smooth two-dimensional manifold for each $\varepsilon<\varepsilon_{0}$. Furthermore, $\mathbf{U}\in\mathcal{W}_{-}^{cu}(\varepsilon)$ with $|P_{-}^{cu}(r_{0})\mathbf{U}(r_{0})|<\delta_{1}$ if and only if
  \begin{align}
      \mathbf{U}(r_{0}) &= \left[d_{1} + \textnormal{O}_{r_{0}}\left(\left[\varepsilon^{2} + |\mathbf{d}|\right]|\mathbf{d}|\right) \right]\mathbf{W}_{1}^{0} + \textnormal{O}_{r_{0}}\left(\left[\varepsilon^{2} + |\mathbf{d}|\right]|\mathbf{d}|\right)\mathbf{W}_{2}^{0} + d_{3}\mathbf{W}_{3}(r_{0}) \nonumber\\
      &\qquad \qquad + \textnormal{O}_{r_{0}}\left(\left[\varepsilon^{2} + |\mathbf{d}|\right]|\mathbf{d}|\right)\mathbf{W}_{4}(r_{0}),\label{McQuighan:Core}
  \end{align}
  for some $\mathbf{d}=(d_{1},d_{3})\in\mathbb{R}^{2}$ with $|\mathbf{d}|<\delta_{1}$, where the right hand side of \eqref{McQuighan:Core} depends smoothly on $(\mathbf{d},\varepsilon)$, and $\textnormal{O}_{r_{0}}(\cdot)$ denotes the standard Landau symbol where the bounding constants may depend on the value of $r_{0}$.\label{Lem:McQ;1}
 \end{Lemma}
 \begin{Proof}
 This follows from Lemma 3.2 and Remark 3.1 in \cite{mcquighan2014oscillons} and involves a standard application of the variation-of-constants formula on a bounded interval $r\in[0,r_{0}]$; see \cite{lloyd2009localized} for a similar example.
 \end{Proof}
 Next, we establish the existence and parametrisation of the far-field manifold $\widetilde{\mathcal{W}}^{s}_{+}(\varepsilon)$, which contains all small-amplitude solutions of \eqref{model:gen}, for sufficiently large $r$, that decay exponentially as $r\to\infty$. To do this, we replace any $\frac{1}{r}$ terms by the unconstrained variable $\sigma$, defined such that \eqref{model:gen} becomes an extended autonomous system of the form,
 \begin{align}
   \begin{pmatrix}
 \mathbf{U} \\\mathbf{V}\\ \sigma
 \end{pmatrix}_{r} &= \begin{pmatrix}
 0 & 1_{2} & 0 \\\mathbf{C} & 0 & 0 \\ 0 & 0 & 0
 \end{pmatrix}\begin{pmatrix}
 \mathbf{U} \\\mathbf{V} \\ \sigma
 \end{pmatrix} + \begin{pmatrix} 0 \\-\sigma \mathbf{V} + \varepsilon^{2}\mathbf{D}\mathbf{U} + \mathbf{f}(\mathbf{U}) \\ -\sigma^{2}
 \end{pmatrix}. \label{model:gen;ext}
 \end{align}
 In order to recover the non-autonomous system \eqref{model:gen}, we would later implement the restriction $\sigma(r_{0}) = r_{0}^{-1}$, thus ensuring $\sigma(r)=r^{-1}$ for all $r\geq r_{0}$. The parametrisation of $\widetilde{\mathcal{W}}_{+}^{s}(\varepsilon)$ relies on the existence of autonomous local centre-stable and centre manifolds for \eqref{model:gen;ext}, as well as the existence of strong stable foliations; this is tackled in great detail in Section 4 of \cite{mcquighan2014oscillons}, which we summarise in the following result,
 \begin{Lemma}
 There exist constants $\varepsilon_{0},\delta_{0},\delta_{1},\rho_{0}>0$ such that the set $\mathcal{W}^{cs}_{+}(\varepsilon)$ of solutions $(U(r),\sigma(r))$ of \eqref{model:gen;ext} for which $\sup_{r\geq\rho_{0}}\left\{\left|U(r)\right|, |\sigma(r)|\right\}<\delta_{0}$ is a smooth four-dimensional manifold for each $0<\varepsilon\leq\varepsilon_{0}$. Furthermore, $(U,\sigma)\in\mathcal{W}^{cs}_{+}(\varepsilon)$ has the expansion
 \begin{align}
     U(r) &= \left[A(r) + \textnormal{O}_{\sigma}\left(|d_{s}|\left[|d_{s}| + |A|+|B|+\varepsilon^{2}\right]\right)\right]\mathbf{W}_{1}^{0}\nonumber\\
     &\qquad + \left[B(r) + \textnormal{O}_{\sigma}\left(|d_{s}|\left[|d_{s}| + |A|+|B|+\varepsilon^{2}\right]\right)\right]\mathbf{W}_{2}^{0}\label{McQuighan:Foliation}\\
     &\qquad + \left[\textnormal{O}_{\sigma}\left(\left[|d_{s}| + |A|+|B|\right]\left[|d_{s}| + |A|+|B|+\varepsilon^{2}\right]\right)\right]\mathbf{W}_{3}(r)\nonumber\\
     &\qquad + \left[d_{s} + \textnormal{O}_{\sigma}\left(\left[ |A|+|B|+\varepsilon^{2}\right]\left[ |A|+|B|+\varepsilon^{2}\right]\right)\right]\mathbf{W}_{4}(r),\nonumber
 \end{align}
 for some $\mathbf{A}=(A,B,d_{s})\in\mathbb{R}^{3}$ with $|\mathbf{A}|<\delta_{1}$ for all $r\geq\rho_{0}$.\label{Lem:McQ;2}
 \end{Lemma}
 \begin{Proof}
 The existence and smoothness of the relevant manifolds and foliations follow from Propositions 4.1 \& 4.2, as well as Lemmas 4.4 \& 4.5 in \cite{mcquighan2014oscillons}. The parametrisation of the local centre-stable manifold is then derived in Section 4.4 of \cite{mcquighan2014oscillons}, resulting in (4.9).
 \end{Proof}
 We note that the functions $A$ and $B$ are defined in \cite{mcquighan2014oscillons} as the projections of a point $(p,\sigma)$ on a local centre manifold $\mathcal{W}^{c}_{+}(\varepsilon)$ onto $\mathbf{W}_{1}^{0}$ and $\mathbf{W}_{2}^{0}$, respectively. In order to parametrise the far-field manifold $\widetilde{\mathcal{W}}_{+}^{s}(\varepsilon)$, we need to determine an expression for $(U,\sigma)\in\mathcal{W}^{cs}_{+}(\varepsilon)$ such that $|U(r)|$ decays exponentially as $r\to\infty$. Assuming that we can determine the conditions such that $A$ and $B$ are exponentially decaying, and so $|U(r)|$ is also exponentially decaying in \eqref{McQuighan:Foliation}, we then evaluate our parametrisation at $\sigma=\frac{1}{r}$ in order to recover solutions $U(r)\in\widetilde{\mathcal{W}}_{+}^{s}(\varepsilon)$ for the non-autonomous system \eqref{model:gen}. Following the parametrisation of the far-field manifold $\widetilde{\mathcal{W}}^{s}_{+}(\varepsilon)$ by exponentially-decaying coordinates $A,B$ on some local centre manifold $\mathcal{W}^{c}_{+}(\varepsilon)$, we now reduce \eqref{model:gen;ext} to a system of equations for the centre coordinates $(A,B,\sigma)$,
 \begin{Lemma}
  Using the coordinates $P^{c}_{+}U(r) = A(r) \mathbf{W}_{1}^{0} + B(r)\mathbf{W}_{2}^{0}$, the vector field \eqref{model:gen;ext} restricted to $\mathcal{W}^{c}_{+}(\varepsilon)$ can be written as 
  \begin{align}
     A_{r} &= B + \mathcal{R}_{A}(A,B,\sigma,\varepsilon),\nonumber\\
     B_{r} &= -\sigma B + c_{0}\left[\varepsilon^{2} A + c^{0}_{2} A^{2}\right]  + \mathcal{R}_{B}(A,B,\sigma,\varepsilon),\label{McQuighan:Normal;eqns}\\
     \sigma_{r} &= -\sigma^{2},\nonumber
 \end{align}
  where $c_{0} = \left(\frac{1 + a}{1 + a^{2}}\right)$ and $c_{2} = \left(1 - a b\right)$. The remainder terms $\mathcal{R}_{A}$ and $\mathcal{R}_{B}$ satisfy
  \begin{align}
      \mathcal{R}_{A} &= \textnormal{O}\left(\left(\varepsilon^{2} + |A|+|B|\right)\left(|A|+|B|\right)\right),\nonumber\\
      \mathcal{R}_{B} &= \textnormal{O}\left(\left(\varepsilon^{2} + |\sigma|^{2}\right)\left(\varepsilon^{2} + |A|\right)|A| + \left(\varepsilon^{2} + |\sigma|^{2} +|A|\right)|\sigma B| + \left(|A|+|B|\right)^{3}\right),\nonumber
  \end{align}
  \label{Lem:McQ;3}
 \end{Lemma}
 \begin{Proof}
 This follows identically to Lemma 4.6 in \cite{mcquighan2014oscillons} except that we no longer have $\mathbb{Z}_{2}$-symmetry, due to the quadratic nonlinearity in \eqref{model:gen}. Instead, we retain the more general normal form seen in \cite[(58) on page 42]{scheel2003radially},
 \begin{align}
     A_{r} &= \widetilde{B},\nonumber\\
     \widetilde{B}_{r} &= -\sigma \widetilde{B} + \gamma_{1}(\varepsilon,\sigma)A + \gamma_{2}(\varepsilon,\sigma)A^{2} + \textnormal{O}\left(\left(\varepsilon^{2} + |\sigma|^{2} +|A|\right)|\sigma \widetilde{B}| + \left(|A|+|\widetilde{B}|\right)^{3}\right),\nonumber\\
     \sigma_{r} &= -\sigma^{2},\nonumber
 \end{align}
 where \cite[Lemma 3.9]{scheel2003radially} implies that
 \begin{align}
     \gamma_{1}(\varepsilon,\sigma) = \varepsilon^{2}\left[c_{0} + \textnormal{O}\left(|\sigma|^{2} + \varepsilon^{2}\right)\right], \qquad \gamma_{2}(\varepsilon,\sigma) = c_{0}\left[c^{0}_{2} + \textnormal{O}\left(|\sigma|^{2} + \varepsilon^{2}\right)\right],\nonumber
 \end{align}
 for some fixed constants $c_{0},c_{2}^{0}$. In order to obtain the normal form \eqref{McQuighan:Normal;eqns}, we introduce a new coordinate $B$ such that $\widetilde{B}=B + \mathcal{R}_{A}(A,B,\sigma,\varepsilon)$. Then, because $\mathcal{R}_{A}$ is higher order in $\widetilde{B}$, no new terms are introduced to leading order and we obtain the leading-order expansion seen in \eqref{McQuighan:Normal;eqns}. Finally, the leading order expansions of each coefficient can be computed via weakly nonlinear analysis seen in \cite{Burke2008Classification}; see Appendix \ref{App:1,2}. We note that these coefficient values can also be found rigorously via normal form transformations, as seen in \cite{scheel2003radially}. 
 \end{Proof}
 In Section 5 of \cite{mcquighan2014oscillons}, localised solutions were found to take the form $A(r) = \varepsilon A_{2}(\varepsilon r)$, where $A_{2}(s)$ satisfies
 \begin{align}
     \left(\partial_{ss} + \frac{1}{s}\partial_{s}\right)A_{2} = A_{2} + c_{3}^{0}A_{2}^{3},\label{McQ:GL}
 \end{align}
 where $c_{3}^{0}$ is the cubic coefficient in the radial normal form, such that the full system \eqref{McQuighan:Syst} has a localised solution if there exists a localised solution to the real non-autonomous Ginzburg-Landau equation \eqref{McQ:GL}. Similarly, we can take an ansatz of the form $A(r) = \varepsilon^{2} A_{2}(\varepsilon r)$, such that the quadratic nonlinearity scales at the same rate as the linear part of the normal form \eqref{McQuighan:Normal;eqns}. Then, we see that localised solutions to \eqref{McQuighan:Normal;eqns} are associated with localised solutions to the real non-autonomous Ginzburg-Landau equation
 \begin{align}
     \left(\partial_{ss} + \frac{1}{s}\partial_{s}\right)A_{2} = c_{0}\left[A_{2} + c_{2}^{0}A_{2}^{2}\right].\label{Case12:GL}
 \end{align}
 We note that \eqref{Case12:GL} can be rescaled such that 
 \begin{align}
     \left(\partial_{\widetilde{s}\widetilde{s}} + \frac{1}{\widetilde{s}}\partial_{\widetilde{s}}\right)\widetilde{A}_{2} = \widetilde{A}_{2} - \widetilde{A}_{2}^{2},\label{Case12:GL;resc}
 \end{align}
 where we have defined $\widetilde{s}:=\sqrt{c_{0}}s$, $\widetilde{A}_{2}:=-c_{2}^{0}A_{2}$. By \cite{Kwong1989}, we know that there exists a unique positive localised solution to \eqref{Case12:GL;resc}; we note that $\textnormal{sgn}\,(A) = -\textnormal{sgn}\,(c_{2}^{0})$, and so we expect to find localised depressed spikes. However there remains some significant technical analysis required to complete this existence proof: including setting up appropriate geometric blow-up charts, proving the existence of a connecting orbit between the two charts related to localised solutions of \eqref{Case12:GL}, and tracking solutions through the various charts. Hence, the full existence proof is beyond the scope of this paper. We will, however, use this as motivation for our numerical investigation in Section \ref{s:Numerics}, where we observe localised depressed spikes bifurcating from the homogeneous bifurcation point, both sub- and super-critically. In the sub-critical region the spikes decay to the bare state as $r\to\infty$, whereas the spikes in the super-critical region decay to the uniform vegetated state as $r\to\infty$. This is all in agreement with our analysis, and so we expect that these localised solutions exist and bifurcate from the homogeneous bifurcation point $P=0$ as previously discussed.

For both cases 1 and 2, we have followed the work of \cite{mcquighan2014oscillons} to produce a (currently incomplete) proof of the existence of localised radial spikes bifurcating from the homogeneous bifurcation at $P=0$. The general equation \eqref{model:gen} is sufficiently similar to \cite{mcquighan2014oscillons} that much of the analysis follows identically up to the centre-manifold reduction in the far-field. We expect that the existence of localised solutions to \eqref{model:gen} is dependent on the existence of localised solutions to the Ginzburg-Landau equation \eqref{Case12:GL}, although this analysis has been left for future work. Finally, we applied the result seen in \cite{Kwong1989} to conclude that, when the uniform vegetated state bifurcates super-critically at $P=0$, localised depressed spikes may emerge in both cases 1 and 2. This is supported by our numerical results in which we find localised depressed spikes bifurcating both sub- and super-critically from the homogeneous instability for the full von Hardenberg model \eqref{vonHard:1}; see Section \ref{s:Numerics}. 

We note that in case 1 solutions are bifurcating from the bare state $\mathcal{B}$, and so localised depressed spikes are always unphysical in our vegetation model. In case 2 solutions are bifurcating from the uniform vegetated state $\mathcal{V}$, and so it is possible that these solutions may be physical if the amplitude of the spike is sufficiently small such that $N(r)$ and $W(r)$ remain non-negative. Our results also suggest that models where the uniform vegetated state bifurcates sub-critically, such as \cite{Gilad2007,Meron2007}, may allow for localised elevated spikes to emerge close to the critical precipitation $p=p_{c}$. However, a closer analysis of the model would first be required in order to substantiate such a prediction.
\subsection{Case 3: Turing Bifurcation}
\label{s:Case3}
\begin{figure}[t!]
    \centering
    \includegraphics[width=\textwidth]{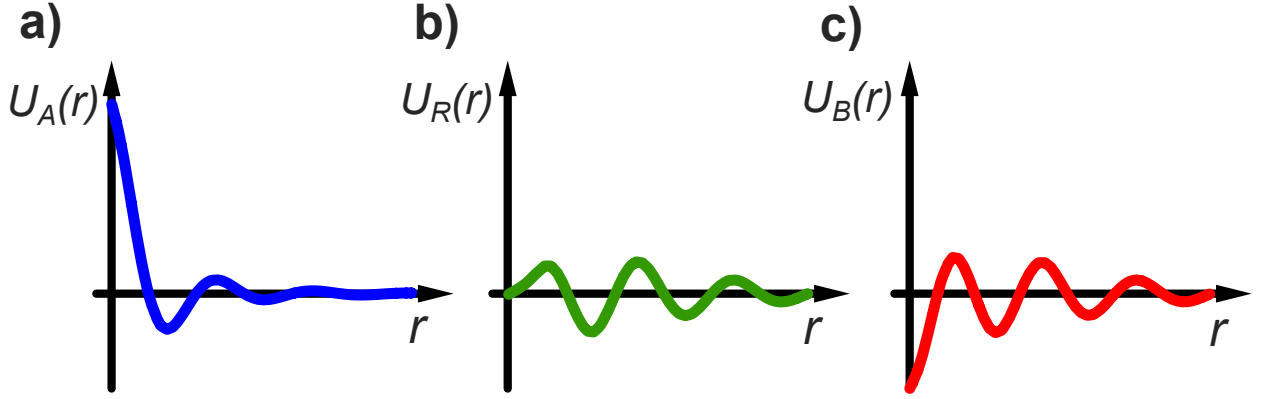}
    \caption{Radial profiles for a) spot A, b) ring, and c) spot B solutions to \eqref{model:3;tilde}. Here, the first element of $U_{A/R/B}(r)$ is plotted for each solution, as defined in a)\eqref{SpotA:U}, b)\eqref{Rings:U}, and c) \eqref{SpotB:U}.}
    \label{fig:RadialPatterns}
\end{figure}

 We recall that case 3 occurs close to a Turing bifurcation, and so we expect that we can relate this case to the study of localised radial patterns in the Swift-Hohenberg equation \cite{lloyd2009localized}. In fact, the system \eqref{model:3;tilde} has been transformed in such a way that the linear part is identical to the Swift-Hohenberg equation seen in \cite{lloyd2009localized}, and so the results of \cite{lloyd2009localized,mccalla2013spots} will follow almost identically. We first note that the linearisation of \eqref{model:3;tilde} at the bifurcation point $\varepsilon=0$
\begin{align}
 \begin{pmatrix}
 U_{1} \\ U_{2}\\ U_{3}\\ U_{4}
 \end{pmatrix}_{r} &= \begin{pmatrix}0 & 0 & 1 & 0\\ 0 & 0 & 0 & 1\\ -1 & 1 & -\frac{1}{\widetilde{r}} & 0\\ 0 & -1 & 0 & -\frac{1}{\widetilde{r}} \end{pmatrix}\begin{pmatrix}
 U_{1} \\ U_{2}\\ U_{3}\\ U_{4}
 \end{pmatrix}, \label{model:3;lin}
\end{align} 
has four linearly independent solutions, as seen in \cite[(2.4)]{lloyd2009localized}, of the form
\begin{align}
    &\mathbf{V}_{1}(r) = \sqrt{2\pi}\left(J_{0}(r),0,-J_{1}(r),0\right)^{\intercal},&\qquad
    &\mathbf{V}_{2}(r) = \sqrt{2\pi}\left(rJ_{1}(r),2J_{0}(r),rJ_{0}(r),-2J_{1}(r)\right)^{\intercal},&\nonumber\\
    &\mathbf{V}_{3}(r) = \sqrt{2\pi}\left(Y_{0}(r),0,-Y_{1}(r),0\right)^{\intercal},&\qquad
    &\mathbf{V}_{4}(r) = \sqrt{2\pi}\left(rY_{1}(r),2Y_{0}(r),rY_{0}(r),-2Y_{1}(r)\right)^{\intercal},&\nonumber
\end{align}
where $J_{\nu}$, $Y_{\nu}$ denote the $\nu^{\textnormal{th}}$ Bessel functions of the first and second kind, respectively. We denote the projection onto the space $\textnormal{span}\{\mathbf{V}_{1}(r),\mathbf{V}_{2}(r)\}$ as $P_{-}^{cu}(r)$, with null space given by $\textnormal{span}\{\mathbf{V}_{3}(r),\mathbf{V}_{4}(r)\}$, for each fixed $r\geq0$. We now conclude the following results from \cite{lloyd2009localized,mccalla2013spots}. The first result is regarding the existence and parametrisation of the core manifold $\mathcal{W}_{-}^{cu}(\varepsilon)$, which contains all small-amplitude solutions of \eqref{model:3;tilde} on the interval $r\in[0,r_{0}]$, for some large fixed $r_{0}>0$, that remain bounded as $r\to0$,
 \begin{Lemma}
  Fix $r_{0}>0$; then, there are constants $\delta_{0},\delta_{1},\varepsilon_{0}>0$ such that the set $\mathcal{W}^{cu}_{-}(\varepsilon)$ of solutions $\mathbf{U}(r)$ of \eqref{model:3;tilde} for which $\sup_{0\leq r\leq r_{0}}|\mathbf{U}|<\delta_{0}$ is a smooth two-dimensional manifold for each $\varepsilon<\varepsilon_{0}$. Furthermore, $\mathbf{U}\in\mathcal{W}_{-}^{cu}(\varepsilon)$ with $|P_{-}^{cu}(r_{0})\mathbf{U}(r_{0})|<\delta_{1}$ if and only if
  \begin{align}
      \mathbf{U}(r_{0}) &= d_{1}\mathbf{V}_{1}(r_{0}) + d_{2}\mathbf{V}_{2}(r_{0}) + \textnormal{O}_{r_{0}}\left(\left[\varepsilon^{2} + |\mathbf{d}|\right]|\mathbf{d}|\right)\mathbf{V}_{3}(r_{0}) \nonumber\\ & \qquad\qquad + \left[\left(\nu + \textnormal{O}\left(r_{0}^{\frac{1}{2}}\right)\right)d_{1}^{2} + \textnormal{O}_{r_{0}}\left(\varepsilon^{2}|\mathbf{d}| + |d_{2}|^{2} + |d_{1}|^{3}\right)\right]\mathbf{V}_{4}(r_{0}),\label{Lloyd:Core}
  \end{align}
  for some $\mathbf{d}=(d_{1},d_{2})\in\mathbb{R}^{2}$ with $|\mathbf{d}|<\delta_{1}$, where $\nu=\omega\sqrt{\frac{\pi}{6}}$, the right hand side of \eqref{Lloyd:Core} depends smoothly on $(\mathbf{d},\varepsilon)$, and $\textnormal{O}_{r_{0}}(\cdot)$ denotes the standard Landau symbol where the bounding constants may depend on the value of $r_{0}$.\label{Lem:Lloyd;1}
 \end{Lemma}
 \begin{Proof}
 This follows from \cite[Lemma 1]{lloyd2009localized} and involves a standard application of the variation-of-constants formula on a bounded interval $r\in[0,r_{0}]$. As seen in \cite{lloyd2009localized}, the quadratic coefficient $\nu$ is found from a Taylor expansion to be
 \begin{align}
    \nu &= \frac{\sqrt{2\pi}}{8}\int_{0}^{r_{0}} r J_{0}(r) \left(\omega\left[\mathbf{V}_{1}(r)\right]_{1}^{2}\right) \,\textnormal{d} r,\nonumber\\
    \intertext{where $\left[\mathbf{X}\right]_{1}$ denotes the first element of a vector $\mathbf{X}$. Hence, we see that}
 \nu&= \frac{\omega\,(2\pi)^{\frac{3}{2}}}{8}\int_{0}^{\infty} r J^{3}_{0}(r)  \,\textnormal{d} r + \textnormal{O}\left(r_{0}^{-\frac{1}{2}}\right)= \omega\sqrt{\frac{\pi}{6}} + \textnormal{O}\left(r_{0}^{-\frac{1}{2}}\right),\nonumber
\end{align}
where we have used \cite[(3) on page 411]{watson1944bessel} to compute the final integral.
 \end{Proof}
 Next, we establish the existence and parametrisation of the far-field manifold $\mathcal{W}^{s}_{+}(\varepsilon)$, which contains all small-amplitude solutions of \eqref{model:3;tilde}, for sufficiently large $r$, that decay exponentially as $r\to\infty$. To do this, we replace any $\frac{1}{r}$ terms by the unconstrained variable $\sigma$, defined such that \eqref{model:3;tilde} becomes an extended autonomous system of the form,
 \begin{align}
 \begin{pmatrix}
 U_{1} \\ U_{2}\\ U_{3}\\ U_{4} \\ \sigma
 \end{pmatrix}_{r} &= \begin{pmatrix}0 & 0 & 1 & 0 & 0\\ 0 & 0 & 0 & 1 & 0 \\ -1 & 1 & 0 & 0 & 0 \\ 0 & -1 & 0 & 0 & 0\\
 0 & 0 & 0 & 0 & 0\end{pmatrix}\begin{pmatrix}
 U_{1} \\ U_{2}\\ U_{3}\\ U_{4}\\ \sigma
 \end{pmatrix} + \begin{pmatrix} 0 \\ 0 \\ -\sigma U_{3} + \left( U_{1} - U_{2}\right)\left(\varepsilon^{2} -U_{1}\right) \\ -\sigma U_{4} -\omega\left( U_{1} - U_{2}\right)\left(\varepsilon^{2} -U_{1}\right) \\-\sigma^{2}
 \end{pmatrix}. \label{model:3;ext}
\end{align} 
 In order to recover the non-autonomous system \eqref{model:3;tilde}, we would later implement the restriction $\sigma(r_{0}) = r_{0}^{-1}$, thus ensuring $\sigma(r)=r^{-1}$ for all $r\geq r_{0}$. The far-field normal form equations are constructed identically to \cite{lloyd2009localized}; we introduce the coordinates $\widetilde{A},\widetilde{B}\in\mathbb{C}$, defined in \cite[(3.5)]{lloyd2009localized} as
 \begin{align}
     \begin{pmatrix}
     \widetilde{A}\\ \widetilde{B}
     \end{pmatrix} &= \frac{1}{4}\begin{pmatrix}
     2U_{1} - \textnormal{i}\left(2U_{3} + U_{4}\right) \\ -U_{4} - \textnormal{i}U_{2}
     \end{pmatrix},\label{coord:AB}
 \end{align}
 and present the following result,
 \begin{Lemma}
 Fix $0<m<\infty$; then, there exists a change of coordinates
 \begin{align}
     \begin{pmatrix}
     A \\ B
     \end{pmatrix} &= \textnormal{e}^{-\textnormal{i}\phi(r)}\left[1 + \mathcal{T}(\sigma)\right]\begin{pmatrix}
     \widetilde{A} \\ \widetilde{B}
     \end{pmatrix} + \textnormal{O}\left(\left(\varepsilon^{2} + |\widetilde{A}| + |\widetilde{B}|\right)\left(|\widetilde{A}| + |\widetilde{B}|\right)\right),\label{Lloyd:Normal;transf}
 \end{align}
 so that the \eqref{model:3;ext}, evaluated at \eqref{coord:AB}, becomes
 \begin{align}
     A_{r} &= -\frac{\sigma}{2} A + B + \mathcal{R}_{A}(A,B,\sigma,\varepsilon),\nonumber\\
     B_{r} &= -\frac{\sigma}{2} B + c_{0}\varepsilon^{2} A + c_{3} |A|^{2}A  + \mathcal{R}_{B}(A,B,\sigma,\varepsilon),\label{Lloyd:Normal;eqns}\\
     \sigma_{r} &= -\sigma^{2}.\nonumber
 \end{align}
 The coordinate change \eqref{Lloyd:Normal;transf} is polynomial in $(A,B,\sigma)$ and smooth in $\varepsilon$, and $\mathcal{T}(\sigma)=\textnormal{O}(\sigma)$ is linear and upper triangular for each $\sigma$, while $\phi(r)$ satisfies
 \begin{align}
     \phi_{r} &= 1 + \textnormal{O}\left(\varepsilon^{2} + |\sigma|^{3} + |A|^{2}\right), \qquad \phi(0)=0.\nonumber
 \end{align}
 The constants $c_{0},c_{3}$ are given by $c_{0}=\frac{\omega}{4}$ and $c_{3}=\frac{\omega(23\omega-30)}{36}$, and the remainder terms satisfy
 \begin{align}
     \mathcal{R}_{A}(A,B,\sigma,\varepsilon) &= \textnormal{O}\left(\sum_{j=0}^{2}\left|A^{j}B^{3-j}\right| + |\sigma|^{3}|A| + |\sigma|^{2}|B| + \left(|A|+|B|\right)^{5} + \varepsilon^{2}|\sigma|^{m}\left(|A|+|B|\right)\right),\nonumber\\
     \mathcal{R}_{B}(A,B,\sigma,\varepsilon) &= \textnormal{O}\left(\sum_{j=0}^{1}\left|A^{j}B^{3-j}\right| + |\sigma|^{3}|B| + \varepsilon^{2}\left(\varepsilon^{2} + |\sigma|^{3} + |A|^{2}\right)|A|+ \left(|A|+|B|\right)^{5} + \varepsilon^{2}|\sigma|^{m}|B|\right).\nonumber
 \end{align}
 \label{Lem:Lloyd;Normal}
 \end{Lemma}
 \begin{Proof}
 This follows identically from \cite[Lemma 2]{lloyd2009localized}, where we compute the leading order expansions of $c_{0},c_{3}$ via weakly nonlinear analysis seen in \cite{Burke2008Classification}; see Appendix \ref{App:3}. Again, we note that these coefficient values can also be found rigorously via normal form transformations, as seen in \cite{scheel2003radially}.
 \end{Proof}
 Localised solutions of \eqref{Lloyd:Normal;eqns} are found to take the form $A(r) = \varepsilon^{\frac{1}{2}} A_{2}\left(\varepsilon^{\frac{1}{2}} r\right)$, where $A_{2}(s)$ satisfies
 \begin{align}
     \left(\partial_{s} + \frac{1}{2s}\right)^{2}A_{2} = c_{0}A_{2} + c_{3} |A_{2}|^{2} A_{2}.\label{Lloyd:GL}
 \end{align}
It was determined in \cite{lloyd2009localized,mccalla2013spots} that there are two types of localised solutions to \eqref{Lloyd:GL} for $c_{0}>0$: one that satisfies the full nonlinear equation \eqref{Lloyd:GL}, and another where $|A_{2}|\ll1$ and so solutions follow the linear flow of \eqref{Lloyd:GL}. The linear flow of \eqref{Lloyd:GL} has a unique localised solution that can be found explicitly and is independent of $c_{3}$; it was proven in \cite{vandenberg2015Rigorous} that the cubic equation \eqref{Lloyd:GL} has a unique positive localised solution $q(s)\textnormal{e}^{\textnormal{i}Y}$, for some $Y\in\mathbb{R}$, if and only if $c_{3}<0$. The real function $q(s)$ has the asymptotic form
\begin{align}
    q(s)=\left\{\begin{array}{cc}
        q_{0}s^{\frac{1}{2}} + \textnormal{O}\left(s^{\frac{3}{2}}\right), & s\to 0 \\
         \left(q_{+} + \textnormal{O}\left(\textnormal{e}^{-\sqrt{c_{0}}s}\right)\right)s^{-\frac{1}{2}}\textnormal{e}^{-\sqrt{c_{0}}s}, & s\to\infty,
    \end{array}\right.\nonumber
\end{align}
where $q_{0}>0$ and $q_{+}\neq0$ are fixed constants. By introducing geometric blow-up coordinates, as seen in \cite{mccalla2013spots}, one can track small-amplitude exponentially decaying solutions back through extended phase spaces by following close the localised solutions of \eqref{Lloyd:GL}. This is entirely identical to the work in \cite{lloyd2009localized,mccalla2013spots}, which allows us to state the following theorems:
\begin{Theorem}[Existence of Spot A]
Fix $\omega>0$. Then, there is some $\varepsilon_{0}>0$ such that \eqref{model:3;tilde} has localised radial solutions $\mathbf{V}_{A}(r):=\left(U_{A}, V_{A}\right)^{\intercal}(r)$, where $U_{A},V_{A}\in\mathbb{R}^{2}$ for all $r\geq0$, of amplitude $\textnormal{O}\left(\varepsilon\right)$ for each $\varepsilon\in(0,\varepsilon_{0})$. In particular, $U_{A}(r)$ has the expansion
\begin{align}
    U_{A}(r) &\sim \frac{\sqrt{3}}{\omega} \varepsilon\begin{pmatrix}
    J_{0}(r) \\ 0 
    \end{pmatrix} + \textnormal{O}(\varepsilon^{2}),\label{SpotA:U}
\end{align}
uniformly on bounded intervals $[0,r_{0}]$ as $\varepsilon\to0$, where  $J_{0}(r)$ is a Bessel function of the first kind of order $0$. \label{Thm:Lloyd;SpotA}
\end{Theorem}
\begin{Theorem}[Existence of Rings]
Fix $\omega>\omega_{*}:=\frac{30}{23}$. Then, there is some $\varepsilon_{0}>0$ such that \eqref{model:3;tilde} has localised radial solutions $\mathbf{V}^{\pm}_{R}(r):=\left(U^{\pm}_{R}, V^{\pm}_{R}\right)^{\intercal}(r)$, where $U^{\pm}_{R},V^{\pm}_{R}\in\mathbb{R}^{2}$ for all $r\geq0$, of amplitude $\textnormal{O}\left(\varepsilon^{\frac{3}{2}}\right)$ for each $\varepsilon\in(0,\varepsilon_{0})$. In particular, $U^{\pm}_{R}(r)$ has the expansion
\begin{align}
    U_{R}^{\pm}(r) &\sim \pm q_{0} \varepsilon^{\frac{3}{2}}\begin{pmatrix} 
    rJ_{1}(r) \\ 2 J_{0}(r)
    \end{pmatrix} + \textnormal{O}(\varepsilon^{2}),\label{Rings:U}
\end{align}
uniformly on bounded intervals $[0,r_{0}]$ as $\varepsilon\to0$, where $q_{0}>0$ is a fixed constant and $(J_{1}(r), J_{0}(r))$ are Bessel functions of the first kind of orders $1$ and $0$, respectively. \label{Thm:Lloyd;Rings}
\end{Theorem}
\begin{Theorem}[Existence of Spot B]
Fix $\omega>\omega_{*}:=\frac{30}{23}$. Then, there is some $\varepsilon_{0}>0$ such that \eqref{model:3;tilde} has a localised radial solution $\mathbf{V}_{B}(r):=\left(U_{B}, V_{B}\right)^{\intercal}(r)$, where $U_{B},V_{B}\in\mathbb{R}^{2}$ for all $r\geq0$, of amplitude $\textnormal{O}\left(\varepsilon^{\frac{3}{4}}\right)$ for each $\varepsilon\in(0,\varepsilon_{0})$. In particular, $U_{B}(r)$ has the expansion
\begin{align}
    U_{B}(r) &\sim -\sqrt{\frac{q_{0} \sqrt{3}}{\omega}}\varepsilon^{\frac{3}{4}}\begin{pmatrix}
    J_{0}(r)\\ 0
    \end{pmatrix} + \textnormal{O}(\varepsilon),\label{SpotB:U}
\end{align}
uniformly on bounded intervals $[0,r_{0}]$ as $\varepsilon\to0$, where $q_{0}>0$ is a fixed constant and $J_{0}(r)$ is a Bessel function of the first kind of order $0$. 
\label{Thm:McCalla;spotB}
\end{Theorem}
Each of these theorems follow identically from the work of \cite{lloyd2009localized,mccalla2013spots}, where the corresponding parameter values are computed in Lemmas \ref{Lem:Lloyd;1} and \ref{Lem:Lloyd;Normal}; see Figure \ref{fig:RadialPatterns} for the profile of each solution. We note that $\omega:=\frac{k^{2} + 1}{k^{2}}>1$; hence, spot A solutions exist for all choices of $\omega$, whereas spot B and ring solutions require that $\omega>\omega^{*}\approx 1.3$ in order to emerge.
\begin{Remark}
We can invert the coordinate transformations defined in \eqref{Defn:P} to obtain expressions for the water and vegetation densities $(N,W)(r)$. Then, in the core region $0\leq r \leq k^{-1} r_{0}$, we find that
\begin{align}
    \textnormal{Spot A:} \quad N_{A}(r) &\sim 
    \frac{k^{4}}{1 - a b}\left[1  + \frac{\sqrt{3}}{\omega} \frac{\varepsilon}{k^{2}}
    J_{0}(k r) + \textnormal{O}(\varepsilon^{2})\right],\nonumber\\
     W_{A}(r) &\sim \frac{k^{4}}{b(1 - a b)}\left[1 + \frac{\sqrt{3}}{\omega} \frac{\varepsilon}{k^{2}} J_{0}(k r)
    + \textnormal{O}(\varepsilon^{2})\right] - \frac{k^{4}}{b}\left[1 - \frac{\sqrt{3}}{\omega} \frac{\varepsilon}{k^{4}} J_{0}(k r)
    + \textnormal{O}(\varepsilon^{2})\right],\nonumber\\
    \textnormal{Spot B:} \quad N_{B}(r) &\sim \frac{k^{4}}{1 - a b}\left[  1 -\sqrt{\frac{q_{0} \sqrt{3}}{\omega}}\frac{\varepsilon^{\frac{3}{4}}}{k^{\frac{3}{2}}}
    J_{0}(k r) + \textnormal{O}(\varepsilon)\right],\nonumber\\
    W_{B}(r) &\sim \frac{k^{4}}{b(1 - a b)}\left[1  -\sqrt{\frac{q_{0} \sqrt{3}}{\omega}}\frac{\varepsilon^{\frac{3}{4}}}{k^{\frac{3}{2}}}J_{0}(k r) + \textnormal{O}(\varepsilon)\right] - \frac{k^{4}}{b}\left[1  +\sqrt{\frac{q_{0} \sqrt{3}}{\omega}}\frac{\varepsilon^{\frac{3}{4}}}{k^{\frac{7}{2}}}J_{0}(k r) + \textnormal{O}(\varepsilon)\right],\nonumber\\
    \textnormal{Rings:} \quad N_{R}^{\pm}(r) &\sim \frac{k^{4}}{1 - a b} \left[1 \pm q_{0} \frac{\varepsilon^{\frac{3}{2}}}{k^{2}}
    r J_{1}(k r) + \textnormal{O}(\varepsilon^{2})\right],\nonumber\\
    W_{R}^{\pm}(r) &\sim \frac{k^{4}}{b(1 - a b)}\left[1 \pm q_{0} \frac{\varepsilon^{\frac{3}{2}}}{k^{2}} r J_{1}(k r) + \textnormal{O}(\varepsilon^{2})\right] -\frac{k^{4}}{b}\left[1 \mp q_{0} \frac{\varepsilon^{\frac{3}{2}}}{k^{5}}\left[k r J_{1}(k r) - 2 J_{0}(k r)\right] + \textnormal{O}(\varepsilon^{2})\right].\nonumber
\end{align}
\end{Remark}
Hence, we have found three classes of localised radial patterns bifurcating from the uniform vegetated state at the Turing point $P=P_{0}$. These are all physical solutions for $\varepsilon$ sufficiently small and have oscillatory tails with an associated wave number $k:=\sqrt{\frac{1+\sqrt{1-K}}{-K}}$, for $K=\frac{1-b}{1-a b}$. If these solutions persist away from the Turing point, these oscillatory tails may be damped by the stronger localisation, resulting in solutions that appear to have monotonic tails even though they bifurcate from a Turing instability. 
\section{Numerical Results for the von Hardenberg Model}
\label{s:Numerics}
We numerically investigate the existence of localised spot A and spike solutions in the full von Hardenberg model \eqref{model:wn}. We employ finite-difference and numerical continuation methods in order to find localised solutions to \eqref{model:wn} and determine their bifurcation structure. While we expect to be able to also find localised rings and spot B solutions, as proven to exist in the reduced model in Theorems \ref{Thm:Lloyd;Rings} \& \ref{Thm:McCalla;spotB}, these solutions typically exist in more restrictive parameter regions than spot A solutions and are harder to find numerically. Hence, while we do expect to be able to find such patterns, we restrict our attention to localised spike and spot A solutions. We also investigate the role of the shading parameter $\rho$ in the emergence of these localised radial patterns and their bifurcation structure. In \cite{Gowda2014}, the default choice of parameters is the following,
\begin{align}
     \gamma = 1.6, \quad \sigma = 1.6, \quad \nu = 0.2, \quad \beta = 3, \quad \delta = 30, \quad \rho = 1.5. \label{param:Gowda}
\end{align}
In this parameter regime, peaks arranged in a hexagonal lattice bifurcate sub-critically from $p=p_{1}$ whereas stripes bifurcate super-critically; see Figure 8 in \cite{Gowda2014}. In \cite{Dawes2016}, localised one-dimensional solutions do not emerge for this parameter region; in order to numerically find such solutions, the authors introduce a second choice of parameters,
\begin{align}
     \{\gamma, \sigma, \nu, \beta, \delta\} \;\textnormal{as defined in \eqref{param:Gowda},} \qquad \rho = 2.5, \label{param:Dawes}
\end{align}
which ensures that stripes bifurcate sub-critically, and so localised one-dimensional solutions emerge from the Turing point $p=p_{1}$. We wish to see how the choice of $\rho$ affects localised radial patterns; we fix the other parameters $\{\gamma, \sigma, \nu, \beta, \delta\}$ as in \eqref{param:Gowda} and present our results for four distinct cases:
\begin{align}
      \mathrm{(i)}\; \rho = 1.5, \qquad   \mathrm{(ii)}\; \rho = 2, \qquad   \mathrm{(iii)}\; \rho = 2.5, \qquad   \mathrm{(iv)}\; \rho = 2.7. \label{param:choices}
\end{align}
These cases include the choices of $\rho$ seen in \cite{Gowda2014} and \cite{Dawes2016} (cases (i) and (iii), respectively), as well as an intermediate value (case (ii)) and an exceeding value (case (iv)) for $\rho$. We note that our analysis holds as long as $a<\frac{1}{b}$, where
\begin{align}
    a := \frac{1}{\beta}\left(\rho - p_{c}\right)p_{c}, \qquad b:= \frac{\beta\gamma}{\left(1 + \sigma p_{c}\right)^{2}}, \qquad p_{c}:= \frac{\nu}{\gamma - \sigma \nu},\nonumber
\end{align}
which, for $\{\gamma, \sigma, \nu, \beta, \delta\}$ fixed as in \eqref{param:Gowda}, is equivalent to the condition $\rho<\frac{\left(1 + \sigma p_{c}\right)^{2}}{\gamma p_{c}} + p_{c} \approx 6.41$. Thus, each case in \eqref{param:choices} is valid for our weakly nonlinear analysis.

\paragraph{Numerical Implementation} We briefly describe our numerical scheme for finding localised solutions to \eqref{model:wn} via central finite-difference methods. We recall that the von Hardenberg model \eqref{model:wn} can be written in the form
\begin{align}
    \mathbf{0} &= \mathbf{F}(U;p):= \Delta_r \begin{pmatrix}
    1 & 0 \\ -\delta\beta & \delta
    \end{pmatrix} \begin{pmatrix}
    n \\ w
    \end{pmatrix} + \begin{pmatrix}
    \left[\left(\frac{\gamma w}{1 + \sigma w}\right) - n - \nu\right]n\\
    -(w-p) + (\rho - w)wn
    \end{pmatrix}\label{vonHard:F}
\end{align}
where $U:=(n,w)^T$. We define our finite radial domain $0\leq r\leq r_{*}:=400$ and impose the following Neumann boundary conditions
\begin{align}
    \lim_{r\to r_{*}} \partial_{r} U = \mathbf{0}, \qquad \textnormal{and} \qquad \lim_{r\to0}\partial_{r} U = \mathbf{0}, \quad \implies \quad \lim_{r\to0}\Delta_r U = 2\lim_{r\to0}\partial_{rr}U.\label{vonHard:BD}
\end{align}
 We are looking for solutions that decay exponentially to some uniform state as the radius $r$ increases, which is why we use homogeneous Neumann boundary conditions at the outer boundary; the inner boundary condition is defined such that \eqref{vonHard:F} is well-defined in the limit as $r\to0$. We discretise our radial domain into $T:=2000$ grid points $\{r_{i}\}_{i=1}^{T}$ and define $\mathbf{U}:=(\mathbf{n},\mathbf{w})\in\mathbb{R}^{2T}$, where $\mathbf{n}_i = n(r_i)$ and $\mathbf{w}_i=w(r_i)$ for each $i\in[1,T]$, such that we can write \eqref{vonHard:F} as
\begin{align}
    \mathbf{0} &= \mathcal{F}(\mathbf{U};p):= \mathcal{D}_r \begin{pmatrix}
    1_{T} & 0_{T} \\ -\delta\beta 1_{T} & \delta1_{T}
    \end{pmatrix} \begin{pmatrix}
    \mathbf{n} \\ \mathbf{w}
    \end{pmatrix} + \begin{pmatrix}
    \left[\gamma \mathbf{w}\circ\left(\mathbf{e}_{T} + \sigma \mathbf{w}\right)^{-1} - \mathbf{n} - \nu\mathbf{e}_{T}\right]\circ\mathbf{n}\\
    -(\mathbf{w}-p\mathbf{e}_{T}) + (\rho\mathbf{e}_{T} - \mathbf{w})\circ\mathbf{w}\circ\mathbf{n}
    \end{pmatrix}.\label{vonHard:F;Disc}
\end{align}
Here, $1_{T}$, $0_{T}$ and $\mathbf{e}_{T}$ are the identity matrix, zero matrix and vector of ones for $\mathbb{R}^{T}$, respectively. Furthermore, the $\circ$ notation indicates a Hadamard (or element-wise) product and $\mathcal{D}_{r}:= \mathbf{D}_{2} + \mathbf{R}\mathbf{D}_{1}\in\mathbb{R}^{2T\times2T}$ approximates the operator $\Delta_r$. The matrices $\mathbf{D}_{1}$ and $\mathbf{D}_{2}$ are the respective central finite-difference matrices for first and second derivatives, $\mathbf{R}:=\textnormal{diag}(1, r_{2}^{-1}, r_{3}^{-1},\dots,r_{T}^{-1})$ approximates $1/r$, and the first and last rows of each matrix are chosen such that the boundary conditions \eqref{vonHard:BD} are satisfied.

We use Matlab's trust-region solver \textit{fsolve} to solve the discretised problem \eqref{vonHard:F;Disc} for a localised initial guess, and then employ path-following codes from \cite{avitabile2020Numerical} in order to continue our solutions in parameter space. At each localised solution $\mathbf{U}^{*}:=(\mathbf{n}^{*},\mathbf{w}^{*})$ of \eqref{vonHard:F;Disc}, with $p=p^{*}$, we calculate the radial stability of our solution by computing the eigenvalues of the Jacobian of the discretised system evaluated at the localised state $(\mathbf{n}^{*},\mathbf{w}^{*})$. That is, we solve the eigenvalue problem $\lambda \mathbf{V} = \delta_{\mathbf{U}}\mathcal{F}(\mathbf{U}^{*};p^{*})\mathbf{V}$ for some $\mathbf{V}\in\mathbb{R}^{2T}$, where the Jacobian $\delta_{\mathbf{U}}\mathcal{F}(\mathbf{U}^{*};p^{*})$ is the following,
\begin{align}
    \mathcal{D}_r \begin{pmatrix}
    1_{T} & 0_{T} \\ -\delta\beta 1_{T} & \delta1_{T}
    \end{pmatrix} + \begin{pmatrix}
    \textnormal{diag}\left(\gamma \mathbf{w}^{*}\circ\left(\mathbf{e}_{T} + \sigma \mathbf{w}^{*}\right)^{-1} - 2\mathbf{n^{*}} - \nu\mathbf{e}_{T}\right) & \textnormal{diag}\left(\gamma\left[\left(\mathbf{e}_{T} + \sigma \mathbf{w}^{*}\right)^{-2}\right]\circ\mathbf{n}^{*}\right)\\
    \textnormal{diag}\left((\rho\mathbf{e}_{T} - \mathbf{w}^{*})\circ\mathbf{w}^{*}\right) & \textnormal{diag}\left((\rho\mathbf{e}_{T} - 2\mathbf{w}^{*})\circ\mathbf{n}^{*}-\mathbf{e}_{T}\right)
    \end{pmatrix}.\label{vonHard:F;Jac}
\end{align}
If we find an eigenvalue $\lambda$ of \eqref{vonHard:F;Jac} such that $\Re(\lambda)>0$, we conclude that the solution $\mathbf{U}^{*}$ is unstable with respect to radial perturbations; likewise, if every eigenvalue $\lambda$ is such that $\Re(\lambda)<0$, we conclude that the solution $\mathbf{U}^{*}$ is stable with respect to radial perturbations. We emphasise that radial stability is a necessary, but not sufficient, condition for solutions to be stable in the two-dimensional plane; for example, a solution might be stable with respect to radial perturbations but unstable with respect to perturbations in the $x$-direction. Hence, from these calculations we can only definitively say when solutions are unstable, not stable.

\begin{figure}[t!]
    \centering
    \includegraphics[width=\textwidth]{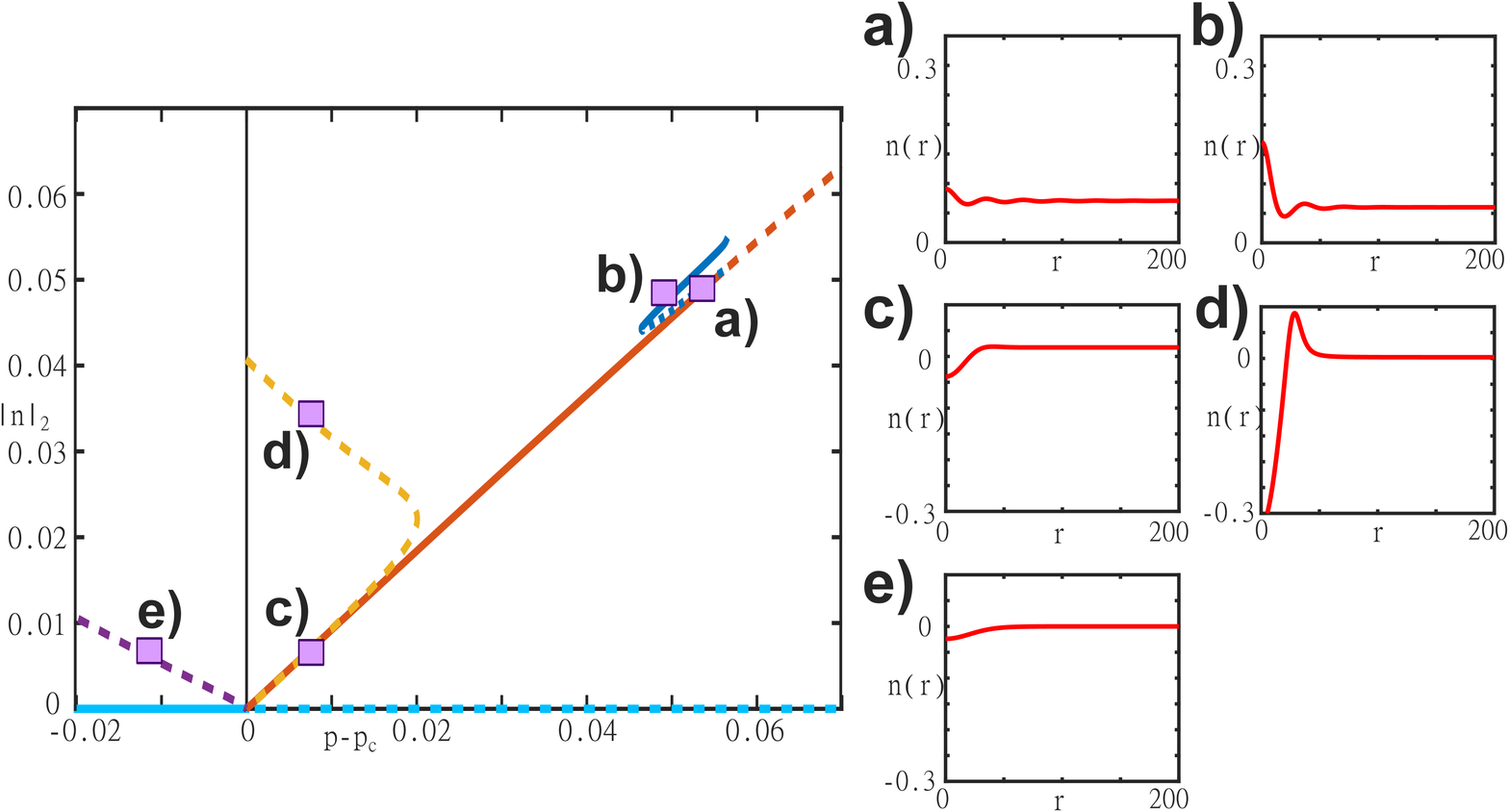}
    \caption{Bifurcation diagram for localised radial solutions to \eqref{model:wn}, when $\rho=1.5$; solid lines indicate stable solutions, and dashed lines unstable solutions. Localised gap solutions bifurcate sub-critically (purple) from the bare state (light blue), as well as super-critically (yellow) from the uniform vegetated state (orange). Localised spot A solutions bifurcate sub-critically (dark blue) from the uniform vegetated state (orange). Panels a) - e) plot the radial profiles for the vegetation density $n(r)$ (red) at the designated points on the bifurcation diagram.}
    \label{fig:rho1.5}
\end{figure}

\begin{figure}[t!]
    \centering
    \includegraphics[width=\textwidth]{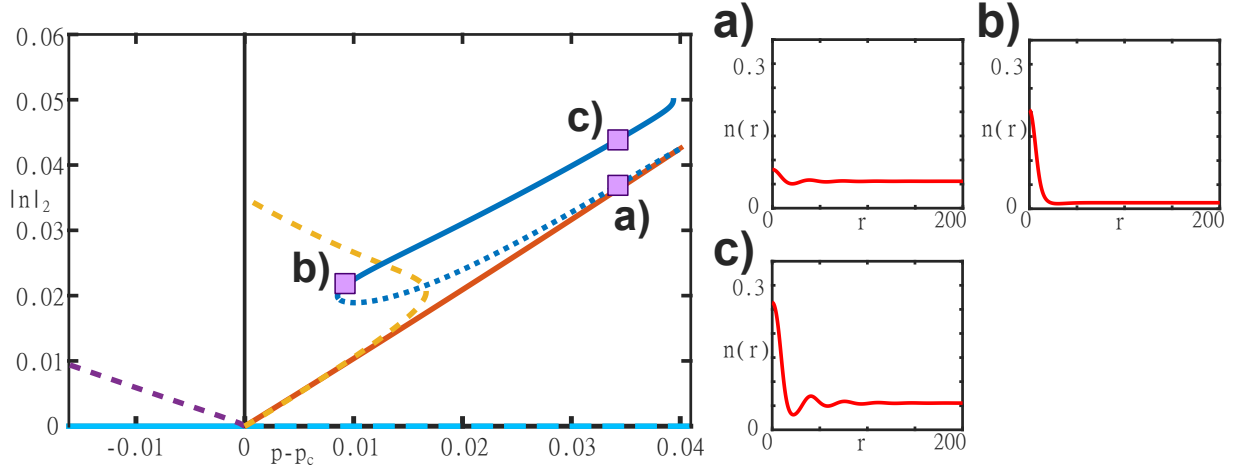}
    \caption{Bifurcation diagram for localised radial solutions to \eqref{model:wn}, when $\rho=2$; solid lines indicate stable solutions, and dashed lines unstable solutions. Localised gap solutions bifurcate sub-critically (purple) from the bare state (light blue), as well as super-critically (yellow) from the uniform vegetated state (orange). Localised spot A solutions bifurcate sub-critically (dark blue) from the uniform vegetated state (orange). Panels a) - c) plot the radial profiles for the vegetation density $n(r)$ (red) at the designated points on the bifurcation diagram.}
    \label{fig:rho2}
\end{figure}

\begin{figure}[t!]
    \centering
    \includegraphics[width=\textwidth]{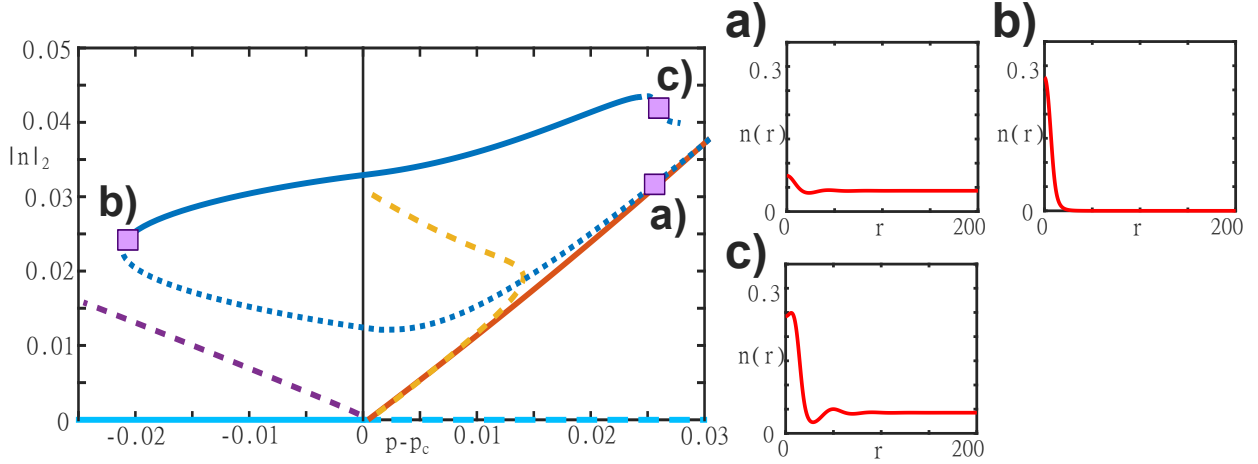}
    \caption{Bifurcation diagram for localised radial solutions to \eqref{model:wn}, when $\rho=2.5$; solid lines indicate stable solutions, and dashed lines unstable solutions. Localised gap solutions bifurcate sub-critically (purple) from the bare state (light blue), as well as super-critically (yellow) from the uniform vegetated state (orange). Localised spot A solutions bifurcate sub-critically (dark blue) from the uniform vegetated state (orange). Panels a) - c) plot the radial profiles for the vegetation density $n(r)$ (red) at the designated points on the bifurcation diagram.}
    \label{fig:rho2.5}
\end{figure}

We note that the size of our radial domain $r_*$ must be chosen sufficiently large in order to find localised solutions. If the domain is too small, then solutions will solve the boundary value problem \eqref{vonHard:F}-\eqref{vonHard:BD} but will not exhibit exponential decay as $r$ increases. Our choice of $r_*=400$ is acceptable in our numerical study, since the observed solutions decay to their uniform state well before reaching the outer boundary. In fact, the solution profiles shown in Figures \ref{fig:rho1.5}-\ref{fig:rho2.7} are only displayed for $r\in[0,200]$; each solution remains flat for the remainder of the radial domain, and so we omit this from the figures. Then, for each localised solution, the profile in the core remains the same for different domain sizes as long as the domain is sufficiently large such that boundary effects are avoided. We now present our numerical findings for cases (i)-(iv).
\paragraph{Case (i)}
The bifurcation diagram for $\rho=1.5$ is plotted in Figure \ref{fig:rho1.5}, where the $L_{2}$-norm of the vegetation density $n(r)$ is plotted against the shifted precipitation $p-p_{c}$. We recall that, in \cite{Dawes2016}, this choice of $\rho$ does not exhibit one-dimensional localised patterns bifurcating from the Turing point $p=p_{1}$. As predicted in our earlier analysis, however, we observe three localised radial patterns: localised gap solutions bifurcating sub- and super-critically from the homogeneous point $p=p_{c}$, respectively, as well as a localised spot A solution bifurcating sub-critically from $p=p_{1}$.

The localised spot A branch is unstable to radial perturbations as it bifurcates from the uniform vegetated state at $p=p_{1}$; as $p$ decreases, the spot undergoes a fold bifurcation where it gains stability with respect to radial perturbations. The changing radial profile of the spot can be seen in panels $a)$ and $b)$ of Figure \ref{fig:rho1.5}. As $p$ re-approaches the critical value $p_1$, the localisation of the spot becomes weaker and its oscillatory tails become more prominent. Eventually the localisation becomes so weak that the oscillating tails reach the outer boundary, and the bifurcation curve of the localised radial spot connects to a domain-covering target pattern; the bifurcation curve for this domain-covering pattern is omitted from Figure \ref{fig:rho1.5}.

The super-critical localised gap solution remains unstable with respect to radial perturbations for all $p>p_{c}$. The solution emerges as a monotonic negative perturbation from the uniform vegetated state at $p=p_{c}$, as seen in Figure \ref{fig:rho1.5} $c)$, and increases in amplitude as $p$ increases. The branch also undergoes a fold bifurcation, although the solution remains unstable to radial perturbations. For the upper part of the solution branch, the localised gap gains a positive peak, as seen in Figure \ref{fig:rho1.5} $d)$. We note that this peak remains much smaller than the depression at $r=0$, however.

The sub-critical localised gap solution bifurcates from the bare state at $p=p_{c}$ and remains unstable with respect to radial perturbations. An example of this solution's radial profile is plotted in Figure \ref{fig:rho1.5} $e)$. As mentioned previously, the sub-critical localised gap is a negative perturbation from the bare state, and so will always be unphysical in the context of our vegetation model. We also note that, for this choice of $\rho$ as well as the remaining cases (ii)-(iv), the super-critical localised gap is a sufficiently large perturbation from the uniform vegetated state that it is also unphysical. As the precipitation increases, this perturbation appears to increase in amplitude faster than the uniform vegetated state, and so we expect the entire branch to be unphysical. Although the localised gap solutions exist in the mathematical model, they are both unphysical for each of our choices of $\rho$; hence, for the remaining cases, we will focus on the behaviour the localised spot A solution.
\paragraph{Case (ii)}
The bifurcation diagram for $\rho=2$ is plotted in Figure \ref{fig:rho2}, with the same measures as in case (i). We note that, in \cite[Figure 15]{Dawes2016}, this choice of $\rho$ may exhibit one-dimensional localised patterns bifurcating from the Turing point $p=p_{1}$, but only for a very small parameter region of $p$. As predicted in our earlier analysis, we observe sub- and super-critical localised gap solutions and a localised spot A solution; see Figure \ref{fig:rho2}.

Similarly to case (i), the localised spot A branch is unstable as it bifurcates and then undergoes a fold bifurcation where it gains stability with respect to radial perturbations. The first difference to note between cases (i) and (ii) is the increased region of existence for the localised spot A solution. By increasing $\rho$ from (i) $1.5$ to (ii) $2$, the length of the spot A branch, in terms of the precipitation $p$, has grown from (i) $0.01$ to (ii) $0.032$. Hence, there is a much larger range of precipitation values for which stable localised spot A solutions exist.

Close to the Turing point $p=p_{1}$, the solution has spatially-oscillating tails with an associated wave number $k_{1}$; see Figure \ref{fig:rho2} $a)$. However, as the solution branch moves further away from $p=p_{1}$, the localisation becomes sufficiently strong that the solution appears to be monotonic; see Figure \ref{fig:rho2} $b)$. We emphasise that oscillations with wavenumber $k_1$ are still present in these solutions, but they have become sufficiently small such that they are not observed qualitatively. This means that a localised spot, typically associated with a Turing bifurcation, can take on an appearance similar to that of a localised spike, typically associated with a homogeneous bifurcation. We note that, following the fold bifurcation, the spot branch returns to $p=p_{1}$, now with a larger amplitude. The localisation of the solution becomes weaker as the precipitation approaches the critical value $p=p_1$, and so spatially-oscillating tails with wave number $k_{1}$ become more prominent again; see Figure \ref{fig:rho2} $c)$. As with case (i), the solution branch for the localised spot terminates when the localisation becomes sufficiently weak that its oscillations reach the outer boundary; here, the solution coincides with a radial target pattern whose bifurcation curve we omit from Figure \ref{fig:rho2}.

\paragraph{Case (iii)}
The bifurcation diagram for $\rho=2.5$ is plotted in Figure \ref{fig:rho2.5}, where the $L_{2}$-norm of the vegetation density $n(r)$ is plotted against the shifted precipitation $p-p_{c}$. We recall that this is the parameter choice in \cite{Dawes2016} such that localised one-dimensional solutions exist and engage in homoclinic snaking; see \cite[Figure 12]{Dawes2016}. Again, we observe sub- and super-critical localised gap solutions and a localised spot A solution; see Figure \ref{fig:rho2.5}.

As with the previous cases, the localised spot A branch is unstable as it bifurcates and then undergoes a fold bifurcation where it gains stability with respect to radial perturbations. Again, compared to the previous cases, we see that there is a larger range of precipitation values for which stable localised spot A solutions exist. A notable change for $\rho=2.5$ is that the localised spot A solution now exists and is stable for some part of the region $p<p_{c}$, where one might expect localised peaks to bifurcate from the homogeneous point $p=p_{c}$ as opposed to the Turing point $p=p_{1}$. As discussed previously, the localisation of the solution is sufficiently strong in this region, such that the spot A may appear to be a monotonic spike even though it has an implicit wave number $k_{1}$. This highlights a potential issue with analysing localised structures experimentally, as their qualitative behaviour may belie the intrinsic properties of their bifurcation.

The behaviour of the localised spot A solution, illustrated in Figure \ref{fig:rho2.5} $a)$ \& $b)$, is identical to case (ii) for most $p$ sufficiently far from the Turing point $p_{1}$. However, for $\rho=2.5$ we note that the upper branch for the localised spot A destabilises close to $p=p_{1}$ in the bifurcation diagram; see Figure \ref{fig:rho2.5}. At this point, the core value decreases such that the maximum value of $n(r)$ is no longer at $r=0$; see Figure \ref{fig:rho2.5} $c)$. This type of behaviour has been seen in \cite[Figure 15]{mcquighan2014oscillons}, in which the radial solution grows outwards and leaves a plateau in the core region. It is unclear whether the same phenomenon would occur in this problem, however, and we leave this question for future study. 

\paragraph{Case (iv)}
Finally, we plot the bifurcation diagram for $\rho=2.7$ in Figure \ref{fig:rho2.7}. We note that this choice of $\rho$ exceeds those used in \cite{Dawes2016}, but we would expect localised one-dimensional solutions to exist and snake within the entire region between the homogeneous point $p_{c}$ and the Turing point $p_{1}$. As with the previous cases, we observe sub- and super-critical localised gap solutions and a localised spot A solution; see Figure \ref{fig:rho2.7}.

We note that the existence region for stable localised spot A solutions is smaller than in case (iii). The stable part of the spot A branch is now almost entirely contained in the region $p<p_{c}$, and so the stable localised spot only appears in the form of a localised monotonic spike; see Figure \ref{fig:rho2.7} $b)$.

The behaviour of the localised spot A solution, illustrated in Figure \ref{fig:rho2.7} $a)$, $b)$ \& $c)$, is identical to case (iii) for most $p$ sufficiently far from the Turing point $p_{1}$. However, for $\rho=2.7$ we note that instability point seen in case (iii) has moved further from the Turing point; see Figure \ref{fig:rho2.7}. At this instability, the core value decreases in a similar way to Figure \ref{fig:rho2.5} $c)$, except with damped spatially-oscillating tails such that they appear monotonic; see Figure \ref{fig:rho2.7} $c)$. Following this unstable branch, the solution decreases in amplitude as it approaches the Turing point and spatially-oscillating tails become more prominent again; see Figure \ref{fig:rho2.7} $d)$. This solution appears similar to a localised ring, suggesting that the spot A branch may connect to other localised radial patterns near the Turing instability, but we leave this question for future study.
\section{Conclusion}
\begin{figure}[t!]
    \centering
    \includegraphics[width=\textwidth]{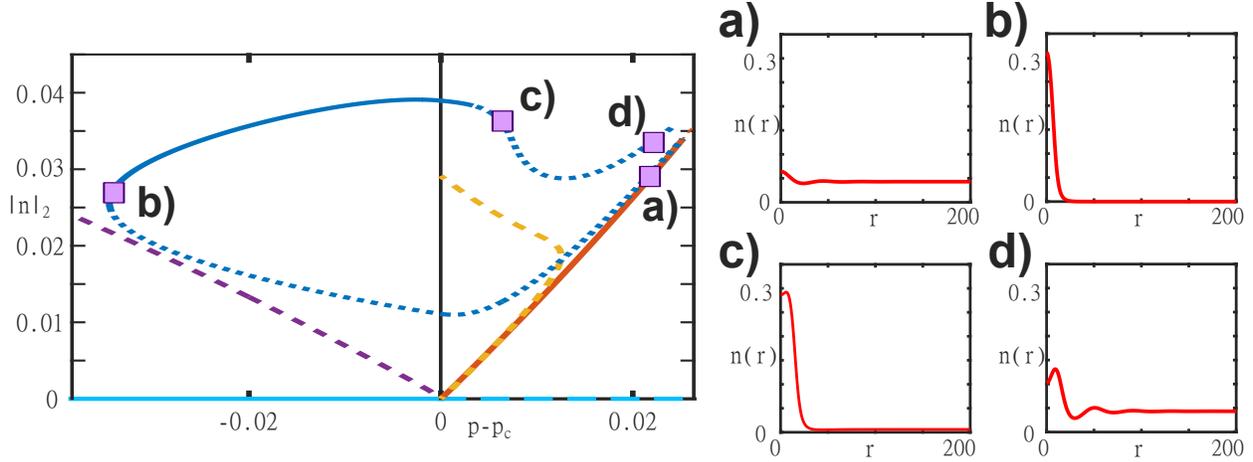}
    \caption{Bifurcation diagram for localised radial solutions to \eqref{model:wn}, when $\rho=2.7$; solid lines indicate stable solutions, and dashed lines unstable solutions. Localised gap solutions bifurcate sub-critically (purple) from the bare state (light blue), as well as super-critically (yellow) from the uniform vegetated state (orange). Localised spot A solutions bifurcate sub-critically (dark blue) from the uniform vegetated state (orange). Panels a) - d) plot the radial profiles for the vegetation density $n(r)$ (red) at the designated points on the bifurcation diagram.}
    \label{fig:rho2.7}
\end{figure}
In this paper, we have investigated the existence of stationary localised radial patterns in the von Hardenberg model for dryland vegetation \cite{vonHardenberg2001}. We began by applying a weakly nonlinear reduction, introduced by \cite{Dawes2016}, to the stationary axisymmetric von Hardenberg model on flat soil \eqref{model:wn}, resulting in the stationary radial reaction-diffusion system \eqref{model:N,W}. By comparing this system to prototypical pattern-forming equations seen in \cite{mcquighan2014oscillons,lloyd2009localized}, we were able to analyse the existence of localised radial patterns close to a homogeneous instability and a Turing instability, respectively.

 In a sub- (super-) critical region close to the homogeneous instability, we presented evidence that localised gap solutions to \eqref{model:N,W} bifurcate from the uniform bare state $\mathcal{B}$ (vegetated state $\mathcal{V}$). We also note that we expect these gap solutions to be unphysical, at least for some range of parameter values, due to the vegetation density becoming negative. The existence proof for these solutions is almost identical to the work presented in \cite{mcquighan2014oscillons}; however, due to a difference in the nonlinearity, further analysis is required to track exponentially decaying solutions in the far-field. This analysis would require tracking solutions through the geometric blow-up charts presented in \cite{mcquighan2014oscillons} for our nonlinear solutions, which is beyond the scope of this paper, and so the proof is left incomplete. Also, in a sub-critical region close to the Turing instability, we proved that three classes of localised radial patterns exist and solve \eqref{model:N,W}. These patterns are small perturbations from the uniform vegetated state $\mathcal{V}$, and are equivalent to the patterns found in the quadratic-cubic Swift-Hohenberg equation \cite{lloyd2009localized,mccalla2013spots}; see Figure \ref{fig:RadialPatterns}.

Following this, we used our analytical results as a guide to numerically find localised radial patterns in the full von Hardenberg model \eqref{model:wn}. For simplicity, we restricted our investigation to the existence and radial stability of localised spot A and gap solutions, which have been predicted by our weakly nonlinear analysis. We used numerical path following routines to continue our solutions to larger amplitudes, and investigated the effect of changing the parameter $\rho$, related to shading effects in \eqref{model:wn}, on the bifurcation structure of the system.

Through this work, we have found that no physically-relevant localised radial patterns bifurcate from the homogeneous instability, at least for a range of appropriate parameter values. Therefore, in this parameter region, any observed localised vegetation patches are predicted to be the result of the uniform vegetated state undergoing a Turing instability. As the precipitation decreases, the oscillatory tails of a localised spot become damped such that they may appear monotonic; in this way, a localised spot can have the appearance of a localised spike far from the Turing instability.

\paragraph{Localised Spots}
We now comment on how our results relate to previous literature regarding localised spots in vegetation models. We first note that the majority of studies of localised spots use models where the uniform vegetated state bifurcates subcritically, and are focused on multi-state solutions, i.e. where there is more than one spot \cite{Meron2007,Gilad2004Engineers,Zelnik2015Gradual}. However, despite these differences, we find that our results remain in qualitative agreement with the wider literature. Localised spots have been observed in the field (see \cite[Fig. 3 \& 7]{Meron2007}, \cite[Fig. 1]{Lejeune2002}) and in many cases appear to be have monotonic tails \cite{Gilad2004Engineers,Zelnik2015Gradual,Meron2007,Lejeune2002}. However, we note that these solutions also tend to exhibit bound states \cite{Meron2007,Lejeune2002}, which suggests that they are localised patterns emerging from a Turing bifurcation \cite{Berrios-Caro2020,Tlidi2020}. This is in harmony with our results, as we only find localised spots bifurcating from the Turing bifurcation and observe that a localised spot can appear monotonic sufficiently far from its associated Turing point; see Figure \ref{fig:rho2} b). 

We also find that localised spots can exist for every choice of shading parameter $\rho$ in \eqref{param:choices} as well as for very low precipitation values, where uniform vegetated states are no longer viable, which is in agreement with results for the Gilad et al. model in \cite{Meron2007}. This is in contrast to localised one-dimensional patterns \cite{Dawes2016}, otherwise known as stripes, which do not emerge for weak shading effects ($\rho<2$) and do not persist for low precipitation values ($p < p_c$). Furthermore, whereas stripes undergo homoclinic snaking, where additional peaks are added to the solution as $p$ oscillates \cite{Dawes2016}, radial spots do not appear to undergo any snaking behaviour and instead terminate at a periodic-pattern branch, as seen in \cite[Fig. 7]{Zelnik2015Gradual} and Figures \ref{fig:rho1.5}-\ref{fig:rho2.7}. 

Another phenomenon observed in the literature is the radial destabilisation of localised spots, resulting in a transition to ring-like structures seen in \cite{Meron2007}. For $\rho=2.5, 2.7$, we observe a radial instability for the localised spot, as seen in Figures \ref{fig:rho2.5} \& \ref{fig:rho2.7}, such that the solution transitions towards a more ring-like structure as the precipitation increases. This is consistent with the behaviour described in \cite[Fig. 10]{Meron2007}, but we leave further study for a follow-up investigation.

\paragraph{Fairy Circles}
Regarding fairy circles, we were unable to find physically-relevant examples of this phenomenon in the von Hardenberg model. However, there may be a number of reasons for this. It may be that the localised depressed spikes bifurcating from the uniform state, pictured as the yellow lines in Figures \ref{fig:rho1.5}-\ref{fig:rho2.7}, can remain strictly positive for a certain choice of parameter values. It could also be that localised spot B solutions, bifurcating from the Turing instability, are responsible for fairy circles. If the spot B solution extends sufficiently far from the Turing instability, much like the spot A solution does, it may begin to resemble a fairy circle more closely. In fact, an analysis of satellite imagery of fairy circles in Namibia and Australia has suggested that such structures exhibit an implicit wavenumber, see \cite[Fig. 3]{Getzin2015Fairy}, which supports our hypothesis that fairy circles are related to the localised spot B solution.

Finally, it may also be possible that the von Hardenberg model is insufficient to find fairy circles, and so a different model may be required. There have been various models used to better study fairy circles; these include the Gilad et al. model \cite{Gilad2007} which, under certain assumptions, can be reduced to a reaction-diffusion system similar to the model studied in this paper \cite{Getzin2015Fairy,Zelnik2015Gradual}. There exist similar reaction-diffusion models for studying fairy circles, as seen in \cite{Zhao2021Fairy}, as well as more complex equations such as the non-local models seen in \cite{Tlidi2020,Berrios-Caro2020,Fernandez-Oto2014interaction}. Some of these models may prove difficult to study analytically, especially those with non-local terms, but any progress could be a major step in developing our understanding of the emergence of fairy circles.

\paragraph{Ecological Implications}
We will briefly discuss some of the ecological implications of this work. We observe that localised spots emerge from the uniform vegetated state $\mathcal{V}$ at a Turing point $p=p_1$; whereas domain-filling patterns mostly exist in the super-critical region $p>p_1$, apart from a possible hysteresis region, localised radial patterns exist in the sub-critical region $p<p_1$. Hence, we expect localised patterns to emerge when the precipitation levels are not sufficient to sustain the domain-covering patterns seen in \cite{Gowda2014}, for example. Furthermore, we expect that localised radial patterns can be found for a much wider range of environments than localised stripes, which are associated with the one-dimensional patterns found in \cite{Dawes2016}. Localised stripes are not found for sufficiently weak shading effects (i.e. for $\rho<2$), suggesting that some species of vegetation may only be able to maintain localised spots and not localised stripes.

For vegetation with weaker shading effects, for example cases (i) and (ii) in \eqref{param:choices}, localised spots appear as a local concentration of vegetation surrounded by a sparser vegetated state. An example of this can be seen in French Guiana, where isolated patches of woody vegetation have been found surrounded by grassland; see \cite[Fig. 1]{Lejeune2002}. For vegetation with stronger shading effects, for example cases (iii) and (iv) in \eqref{param:choices}, we instead find that localised spots persist as the precipitation decreases until they appear as isolated patches of vegetation surrounded by bare soil; examples of this have been found in the Negev region of Israel, such as in \cite[Fig. 3 \& 7]{Meron2007}. Hence, changes to the bifurcation structure caused by the value of the shading parameter $\rho$ can result in the formation of qualitatively different localised spot solutions in the field.

The saddle node of spot A solutions, as seen in Figures \ref{fig:rho1.5}-\ref{fig:rho2.7}, defines a tipping point in climate dynamics; when the precipitation decreases beyond this point, solutions can undergo a transition to the stable uniform state. This transition is hysteretic because the uniform state remains stable as the precipitation returns to its original value. For plant species with sufficiently strong shading effects the spot A saddle node point extends into a region of low precipitation such that the uniform vegetated state is no longer viable; then, this tipping point represents an ecological collapse to the bare state $\mathcal{B}$. This transition is sometimes referred to as a `catastrophic regime shift' \cite{Rietkerk2004catastrophe}, and provides an example of desertification induced by a change in environment. 

\paragraph{Further Study}
A possible extension of this work would be to use the singular perturbation methods seen in \cite{Jaibi2020,Carter2018} to attempt to prove the existence of stationary localised radial patterns in the full von Hardenberg model \eqref{model:wn}. There is also an interesting question regarding the reintroduction of a sloping term into the original model. Instead of setting $c=0$ in \eqref{vonHard:1}, what if we set $c=\hat{c}\delta^{\frac{1}{2}}$ such that there is some small symmetry-breaking that favours stripes over spots? Are we still able to find localised radial solutions for sufficiently small $\hat{c}$ in the reduced model \eqref{model:N,W} and, if so, can we say anything about their stability? There is also the question of whether we can study the existence of more complicated localised patterns, such as solutions with dihedral symmetry \cite{HillPrepPatch} or labyrinthine patterns \cite{Clerc2021Labyrinthine}, in the von Hardenberg model. These are just a couple possibilities for future work; there has also been recent progress on centre manifold reduction techniques for nonlocal equations \cite{faye2018Center}, which we may be able to adapt for radial solutions to nonlocal vegetation models seen in \cite{Meron2007,Tlidi2008}, for example.

\begin{Acknowledgment}
The author would like to thank the anonymous referees for their insightful comments and suggestions, as well as David Lloyd and Matthew Turner for useful discussions on the subject and helpful comments on the manuscript.
\end{Acknowledgment} 
\begin{appendix}
\renewcommand{\theequation}{\Alph{section}.\arabic{equation}}
\section*{Appendix}
\section{Weakly Nonlinear Analysis}
\label{s:App;WeaklyNonlin}
For each case, the radial problem can be written as a generic $2$-dimensional problem,
\begin{align}
    \mathscr{L}\begin{pmatrix}
    u \\ v
    \end{pmatrix} + \mathscr{N}(u,v) &= 0,\label{Weak:Nonlin;gen}
\end{align}
where $\mathscr{L}$ is a linear differential operator and $\mathscr{N}$ contains all the nonlinear terms of the problem such that $\mathscr{N}(0,0)=\mathbf{0}$ and $\delta\mathscr{N}(0,0)=\mathbf{0}$, where $\delta\mathscr{N}$ is the Jacobian of $\mathscr{N}$. Using the work of Scheel \cite{scheel2003radially}, we obtained radial normal forms \eqref{McQuighan:Normal;eqns}, \eqref{Lloyd:Normal;eqns} for the far-field behaviour in cases 1,2 and 3. We now wish to perform weakly nonlinear analysis, as seen in \cite{burke2007homoclinic}, in order to determine certain coefficients for each of the normal forms. In particular, for each case we will derive a nonlinear non-autonomous Ginzburg-Landau equation that determines the behaviour of small-amplitude solutions for large values of $r$.
\subsection{Cases 1 and 2}\label{App:1,2}
For cases 1 and 2, \eqref{Weak:Nonlin;gen} has the following form,
\begin{align}
    &\mathscr{L}= \left(\partial_{rr} + \frac{1}{r}\partial_{r}\right)1_{2} - \mathbf{C} - \varepsilon^{2}\mathbf{D}_{i},& \qquad &\mathscr{N}(u,v)=\begin{pmatrix}
     \left(bv - u\right)u\\ \left(bv - u\right)u 
    \end{pmatrix},&\nonumber
\end{align}    
where,
\begin{align}
    &\mathbf{C}=\begin{pmatrix}
    0 & 0 \\
    - a & 1
    \end{pmatrix},& \qquad &\mathbf{D}_{1}=\begin{pmatrix}
     1 & 0 \\ 1 & 0
    \end{pmatrix},& \qquad &\mathbf{D}_{2}=\begin{pmatrix}
     \frac{1}{1-ab} & -\frac{b}{1-ab} \\ \frac{1}{1-ab} & -\frac{b}{1-ab}
    \end{pmatrix},&\nonumber
\end{align}
respectively. From the normal form for this bifurcation, we expect to find a quadratic Ginzburg-Landau equation that defines the behaviour of solutions in the far-field, after some rescaling. Hence, we define $R=\varepsilon r$ and
\begin{align}
    \begin{pmatrix}
    u(r) \\ v(r)
    \end{pmatrix} &= \varepsilon^{2}\begin{pmatrix}
    u_{1}(R) \\ v_{1}(R)
    \end{pmatrix} + \varepsilon^{4}\begin{pmatrix}
    u_{2}(R) \\ v_{2}(R)
    \end{pmatrix} + \dots\nonumber
\end{align}
Then, we can write $\mathscr{L} = -\mathbf{C} + \varepsilon^{2}\mathscr{L}_{1}$, $\mathscr{N}(u,v) = \varepsilon^{4}\mathscr{N}_{1}(u_{1},v_{1}) + \textnormal{O}(\varepsilon^{6})$, where $\mathscr{L}_{1} = \left(\partial_{RR} + \frac{1}{R}\partial_{R}\right)1_{2} - \mathbf{D}_{i}$, and so \eqref{Weak:Nonlin;gen} can be expressed as
\begin{align}
    &\textnormal{O}(\varepsilon^{2}):& \qquad & \mathbf{C} \begin{pmatrix}
    u_{1}\\v_{1}
    \end{pmatrix} = \mathbf{0}, &\nonumber\\
    &\textnormal{O}(\varepsilon^{4}):& \qquad & \mathbf{C} \begin{pmatrix}
    u_{2}\\v_{2}
    \end{pmatrix} = \left[\mathscr{L}_{1} \begin{pmatrix}
    u_{1}\\v_{1}
    \end{pmatrix} + \mathscr{N}_{1}(u_{1},v_{1})\right]. &\nonumber
\end{align}
The first equation implies that $(u_{1},v_{1})^{\intercal}$ must take the form of an eigenvector of $\mathbf{C}$ with eigenvalue $0$; hence, we define $(u_{1},v_{1})^{\intercal}:= A(R) \widetilde{U}_{0}$, where $\widetilde{U}_{0}$ is defined in \eqref{Defn:U0}. Then, we choose $(u_{2},v_{2})^{\intercal} = B(R)\widetilde{U}_{0}$, such that the second equation becomes
\begin{align}
    \mathscr{L}_{1}\widetilde{U}_{0} A(R) + N(A)\begin{pmatrix}
    1 \\ 1
    \end{pmatrix} &= \mathbf{0},\nonumber
\end{align}
for the function $N(A) = -\left(1 - a b\right) A^{2}(R)$. In order to solve this equation, we define the projection $\mathscr{P}$, where
\begin{align}
    \mathscr{P} = \begin{pmatrix}
    \frac{1}{1+a^{2}} & \frac{a}{1+a^{2}} \\ 
    \frac{a}{1+a^{2}} & \frac{a^{2}}{1+a^{2}} 
    \end{pmatrix}.\nonumber
\end{align}
Then, for both cases 1 and 2, pre-multiplying by $\mathscr{P}$ results in the following equation for $A(R)$,
\begin{align}
    \left(\partial_{RR} + \frac{1}{R}\partial_{R}\right)A &= \left(\frac{1+a}{1+a^{2}}\right)\left[A + \left(1 - a b\right) A^{2} \right].\label{model12:GL}
\end{align}
Hence, we have found that the coefficients $c_{0}$ and $c_{2}$ in \eqref{McQuighan:Normal;eqns} are defined such that $c_{0} = \left(\frac{1+a}{1+a^{2}}\right)$ and $c_{2} = \left(1-a b\right)$. 
\subsection{Case 3}\label{App:3}
For case 3, we write the explicit form of \eqref{Weak:Nonlin;gen} as
\begin{align}
    \mathscr{L}&=\begin{pmatrix}
    \partial_{rr} + \frac{1}{r}\partial_{r} + 1 -  \varepsilon^{2} & \varepsilon^{2} - 1 \\
    \omega\varepsilon^{2} & \partial_{rr} + \frac{1}{r}\partial_{r} + 1 - \omega\varepsilon^{2}
    \end{pmatrix}, \qquad \mathscr{N}(u,v)=\begin{pmatrix}
     \left(u - v\right)u\\ -\omega\left(u - v\right)u 
    \end{pmatrix}.\nonumber
\end{align}
From the normal form for this bifurcation, we expect to find a cubic Ginzburg-Landau equation that defines the far-field behaviour, and so we define $R=\varepsilon r$ and
\begin{align}
    \begin{pmatrix}
    u(r) \\ v(r)
    \end{pmatrix} &= \varepsilon\begin{pmatrix}
    u_{1}(r,R) \\ v_{1}(r,R)
    \end{pmatrix} + \varepsilon^{2}\begin{pmatrix}
    u_{2}(r,R) \\ v_{2}(r,R)
    \end{pmatrix} + \varepsilon^{3}\begin{pmatrix}
    u_{3}(r,R) \\ v_{3}(r,R)
    \end{pmatrix} + \dots\nonumber
\end{align}
Then, we can write $\mathscr{L} = \mathscr{L}_{0} + \varepsilon\mathscr{L}_{1} + \varepsilon^{2}\mathscr{L}_{2}$, $\mathscr{N}(u,v) = \varepsilon^{2}\mathscr{N}_{1}(u_{1},v_{1}) + \varepsilon^{3}\mathscr{N}_{2}(u_{1},v_{1},u_{2},v_{2}) + \textnormal{O}(\varepsilon^{4})$, and so \eqref{Weak:Nonlin;gen} can be expressed as
\begin{align}
    &\textnormal{O}(\varepsilon):& \qquad & \mathscr{L}_{0} \begin{pmatrix}
    u_{1}\\v_{1}
    \end{pmatrix} = \mathbf{0}, &\label{App:3,1}\\
    &\textnormal{O}(\varepsilon^{2}):& \qquad & \mathscr{L}_{0} \begin{pmatrix}
    u_{2}\\v_{2}
    \end{pmatrix} = - \left[\mathscr{L}_{1} \begin{pmatrix}
    u_{1}\\v_{1}
    \end{pmatrix} + \mathscr{N}_{1}(u_{1},v_{1})\right], &\label{App:3,2}\\
    &\textnormal{O}(\varepsilon^{3}):& \qquad & \mathscr{L}_{0} \begin{pmatrix}
    u_{3}\\v_{3}
    \end{pmatrix} = - \left[\mathscr{L}_{1} \begin{pmatrix}
    u_{2}\\v_{2}
    \end{pmatrix} + \mathscr{L}_{2} \begin{pmatrix}
    u_{1}\\v_{1}
    \end{pmatrix} + \mathscr{N}_{2}(u_{1},v_{1},u_{1},v_{2})\right], &\label{App:3,3}
\end{align}
where,
\begin{align}
    &\mathscr{L}_{0} = \begin{pmatrix}
    \partial_{rr} + 1 & -1 \\ 0 & \partial_{rr} + 1
    \end{pmatrix},& \nonumber\\
    &\mathscr{L}_{1} = \begin{pmatrix}
    2\left(\partial_{R} + \frac{1}{2R}\right)\partial_{r} & 0 \\ 0 & 2\left(\partial_{R} + \frac{1}{2R}\right)\partial_{r}
    \end{pmatrix},& \qquad &\mathscr{N}_{1} = \begin{pmatrix}
    (u_{1}-v_{1})u_{1} \\ -\omega(u_{1}-v_{1})u_{1}
    \end{pmatrix},&\nonumber\\
    &\mathscr{L}_{2} = \begin{pmatrix}
    \partial_{RR} + \frac{1}{R}\partial_{R} - 1 & 1 \\ \omega & \partial_{RR} + \frac{1}{R}\partial_{R} - \omega
    \end{pmatrix},& \qquad &\mathscr{N}_{2} = \begin{pmatrix}
    \left[(u_{2}-v_{2})u_{1} + (u_{1}-v_{1})u_{2}\right] \\ -\omega\left[(u_{2}-v_{2})u_{1} + (u_{1}-v_{1})u_{2}\right]
    \end{pmatrix}.& \nonumber
\end{align}
Solving \eqref{App:3,1}, we take 
\begin{align}
    u_{1}(r,R) &= A_{1}(R) \textnormal{e}^{\textnormal{i}r} + \overline{A}_{1}(R) \textnormal{e}^{-\textnormal{i}r}, \qquad v_{1}(r,R) = 0,\nonumber
\end{align}
where overbars denote complex conjugation. Then, \eqref{App:3,2} becomes
\begin{align}
    \begin{pmatrix}
    \left(\partial_{rr} + 1\right)u_{2} - v_{2} \\ \left(\partial_{rr} + 1\right)v_{2}
    \end{pmatrix} &= \begin{pmatrix}
    -2\left(\partial_{R} + \frac{1}{2R}\right)\partial_{r}u_{1} - u^{2}_{1} \\ \omega u^{2}_{1}
    \end{pmatrix}.\nonumber
\end{align}
We take the following ansatz,
\begin{align}
    v_{2}(r,R) &= O_{v,2}(R) + A_{v,2}(R)\textnormal{e}^{\textnormal{i}r} + B_{v,2}(R)\textnormal{e}^{2\textnormal{i}r} + \overline{A}_{v,2}(R)\textnormal{e}^{-\textnormal{i}r} + \overline{B}_{v,2}(R)\textnormal{e}^{-2\textnormal{i}r},\nonumber
\end{align}
such that
\begin{align}
    O_{v,2} - 3 B_{v,2}\textnormal{e}^{2\textnormal{i}r} - 3 \overline{B}_{v,2}\textnormal{e}^{-2\textnormal{i}r} &= \omega \left[A_{1}^{2}\textnormal{e}^{2\textnormal{i}r} + 2|A_{1}|^{2} + \overline{A}_{1}^{2}\textnormal{e}^{-2\textnormal{i}r}\right].\nonumber
\end{align}
Then, equating for each $\textnormal{e}^{a\textnormal{i}r}$, for $a=-2,0,2$, we determine that
\begin{align}
    v_{2}(r,R) &= 2\omega|A_{1}|^{2} + A_{v,2}(R)\textnormal{e}^{\textnormal{i}r} - \frac{\omega}{3} A_{1}^{2}\textnormal{e}^{2\textnormal{i}r} + \overline{A}_{v,2}(R)\textnormal{e}^{-\textnormal{i}r} - \frac{\omega}{3} \overline{A}_{1}^{2}\textnormal{e}^{-2\textnormal{i}r}.\nonumber
\end{align}
Similarly for $u_{2}$, we define
\begin{align}
    u_{2}(r,R) &= O_{u,2}(R) + A_{u,2}(R)\textnormal{e}^{\textnormal{i}r} + B_{u,2}(R)\textnormal{e}^{2\textnormal{i}r} + \overline{A}_{u,2}(R)\textnormal{e}^{-\textnormal{i}r} + \overline{B}_{u,2}(R)\textnormal{e}^{-2\textnormal{i}r},\nonumber
\end{align}
and so
\begin{align}
    O_{u,2} - 3 B_{u,2}\textnormal{e}^{2\textnormal{i}r} - 3 \overline{B}_{u,2}\textnormal{e}^{-2\textnormal{i}r} &= -2\textnormal{i}\left(\partial_{R} + \frac{1}{2R}\right)\left[A_{1}(R) \textnormal{e}^{\textnormal{i}r} - \overline{A}_{1}(R) \textnormal{e}^{-\textnormal{i}r}\right]-\left[A_{1}^{2}\textnormal{e}^{2\textnormal{i}r} + 2|A_{1}|^{2} + \overline{A}_{1}^{2}\textnormal{e}^{-2\textnormal{i}r}\right]\nonumber\\
    &\qquad \qquad + 2\omega|A_{1}|^{2} + A_{v,2}(R)\textnormal{e}^{\textnormal{i}r} - \frac{\omega}{3} A_{1}^{2}\textnormal{e}^{2\textnormal{i}r} + \overline{A}_{v,2}(R)\textnormal{e}^{-\textnormal{i}r} - \frac{\omega}{3} \overline{A}_{1}^{2}\textnormal{e}^{-2\textnormal{i}r}.\nonumber
\end{align}
Hence, we find that $A_{v,2}(R)=2\textnormal{i}\left(\partial_{R} + \frac{1}{2R}\right)A_{1}(R)$ and
\begin{align}
    u_{2}(r,R) &= 2(\omega-1)|A_{1}|^{2} + A_{2}(R)\textnormal{e}^{\textnormal{i}r} + \left(\frac{\omega+3}{9}\right)A_{1}^{2}\textnormal{e}^{2\textnormal{i}r} + \overline{A}_{2}(R)\textnormal{e}^{-\textnormal{i}r} + \left(\frac{\omega+3}{9}\right)\overline{A}_{1}^{2}\textnormal{e}^{-2\textnormal{i}r},\nonumber
\end{align}
which we now need to substitute into \eqref{App:3,3},
\begin{align}
    \begin{pmatrix}
    \left(\partial_{rr} + 1\right)u_{3} - v_{3} \\ \left(\partial_{rr} + 1\right)v_{3}
    \end{pmatrix} &= -\left[2\left(\partial_{R} + \frac{1}{2R}\right)\partial_{r}\begin{pmatrix}
    u_{2} \\ v_{2}
    \end{pmatrix} + \begin{pmatrix}
    \left(\partial_{RR} + \frac{1}{R}\partial_{R}-1\right)u_{1} \\ \omega u_{1}
    \end{pmatrix} + \begin{pmatrix}
    (2u_{2}-v_{2})u_{1}\\ -\omega(2u_{2}-v_{2})u_{1}
    \end{pmatrix}\right].\nonumber
\end{align}
We take the second equation, introduce a similar ansatz to before, and compute the $\textnormal{e}^{\textnormal{i}r}$ terms; this results in the following equation,
\begin{align}
    \left(\partial_{R} + \frac{1}{2R}\right)^{2}A_{1} &= \frac{\omega}{4} A_{1} - \left(\frac{\omega\left(23\,\omega-30\right)}{36} \right)|A_{1}|^{2}A_{1}.\label{model3:GL}
\end{align}
Hence, we find that the coefficients $c_{0}$ and $c_{3}$ in \eqref{Lloyd:Normal;eqns} are defined such that $c_{0}=\frac{\omega}{4}$ and $c_{3} = -\frac{\omega(23\omega-30)}{36}$.
\end{appendix}

\bibliographystyle{abbrv}
\bibliography{Paper.bbl}
\end{document}